\documentclass[10pt]{article}

\usepackage{a4wide}
\usepackage{amssymb}
\usepackage{amsfonts}
\usepackage{amsmath}
\input xy
\xyoption{arrow} \xyoption{matrix}

\date{}

\newtheorem{proposition}{Proposition}[section]
\newtheorem{theorem}[proposition]{Theorem}
\newtheorem{lemma}[proposition]{Lemma}

\newtheorem{corollary}[proposition]{Corollary}

\def\der{\partial }

\def\nFM0{{\nu }_{F,M_0}}
\def\nFN0{{\nu }_{F,N_0}}
\def\nGN0{{\nu }_{G,N_0}}

\def\N0{ {\bf N}_0 }

\def\t{\otimes}
\def\g{\gamma}
\def\v{\varphi}
\def\ra{\rightarrow}

\def\Xpm{X^{\pm }}

\def\s{\sigma}
\def\Z{\mathbb{Z}}

\def\l1{{\lambda}_1}

\def\a{\alpha}
\def\a0{ {\alpha }_0}
\def\a1{ {\alpha }_1}

\def\l{\lambda}
\def\o{\omega}

\def\nFGM0{{\nu }_{F,G,M_0}}

\def\nFN0{{\nu}_{F,N_0}}

\def\sm{{\sigma}^m}

\def\sm1{{\sigma}^{-1}}

\def\smtp1{{\sigma}^{-t+1}}

\def\o{\omega }
\def\S1{S^{-1}}

\def\Xpm1{X^{\pm 1}_1}

\def\sPM1{{\sigma }^{\pm 1}}
\def\sMP1{{\sigma }^{\mp 1 }}

\def\d{\delta}

\def\di{{\rm d.ind}}

\def\L{\Lambda}
\def\O{\Omega}

\def\CD{{\cal D}}

\def\Ytm1{Y^{t-1}}
\def\Yim1{Y^{i-1}}

\def\CM{{\cal M}}

\def\CF{{\cal F}}
\def\CG{{\cal G}}
\def\CH{{\cal H}}

\def\Aut{{\rm Aut}}

\def\Der{{\rm Der }}
\def\ad{{\rm ad }}
\def\dim{{\rm dim }}

\def\ker{ {\rm ker } }

\def\CJ{ {\cal J}}

\def\SL2Z{ {\rm SL}_2({\bf Z}) }

\def\th{ \theta }

\def\Gp1{ G^{1 , 1 } }
\def\P11{ P^{-1 , 1 } }
\def\Pp1{ P^{1 , 1 } }

\def\th{\theta}

\def\nCLsr{{}^\nu\kern-2pt {\cal L}^{\sigma , \rho  }}
\def\nP{{}^\nu \kern-2pt P}
\def\nL{{}^\nu\kern-2pt L}
\def\nLL{{}^\nu\kern-2pt \Lambda}
\def\nPsr{{}^\nu\kern-2pt P^{\sigma , \rho  }}
\def\nLsr{{}^\nu\kern-2pt L^{\sigma , \rho  }}
\def\nuCL{{}^\nu\kern-2pt  {\cal L}}
\def\nCLsr{{}^\nu\kern-2pt {\cal L}^{\sigma , \rho  }}
\def\nCL1m{{}^\nu\kern-2pt {\cal L}^{-1 , 1  }}
\def\x1nu{x^\frac{1}{\nu}}
\def\xm1nu{x^{-\frac{1}{\nu}}}

\def\ra{\rightarrow }

\def\CB{{\cal B}}

\def\CI{{\cal I}}

\def\coker{{\rm coker}}

\def\CC{ {\cal C}}

\def\CH{ {\cal H}}
\def\CP{ {\cal P}}

\def\nAM0{{\nu }_{{\cal A},M_0}}
\def\nAN0{{\nu }_{{\cal A},N_0}}

\def\End{ {\rm End }}
\def\Der{ {\rm Der }}
\def\CJ{ {\cal J }}

\def\CP{ {\cal P }}

\def\ad{ {\rm ad }}

\def\ga{\mathfrak{a}}
\def\gb{\mathfrak{b}}

\def\SL{{\rm SL}}

\def\di!{\frac{\der^i}{i!}}
\def\dik!{\frac{\der^k_i}{k!}}

\def\ord{{\rm ord}}

\def\Nn{\mathbb{N}^n}

\def\gl{\mathfrak{l}}

\def\N{\mathbb{N}}

\def\0{\overline{0}}
\def\1{\overline{1}}

\def\Ln1{\L_{n,\overline{1}}}

\def\a1{a_{\overline{1}}}

\def\S{\Sigma}

\def\bCJ{\overline{\CJ}}

\def\vn1{\overrightarrow{n-1}}

\def\gl{{\rm gl}}

\def\bu{\overline{u}}

\def\Inn{{\rm Inn}}

\def\mJ{\mathbb{J}}
\def\mI{\mathbb{I}}

\def\ann{{\rm ann}}

\def\Cen{{\rm Cen}}

\def\ind{{\rm ind}}

\def\K1{{\rm K}_1}

\def\hmI1{\widehat{\mI_1}}
\def\tmI1{\widetilde{\mI_1}}
\def\tmJ1{\widetilde{\mJ_1}}
\def\hB1{\widehat{B_1}}
\def\hCB1{\widehat{\CB_1}}

\def\mW{\mathrm{W}}

\def\ggu{\mathfrak{u}}
\def\Xai{X_{\alpha , i}}
\def\Xbj{X_{\beta , j}}
\def\udim{{\rm u.dim}}
\def\cdim{{\rm c.dim}}
\def\UT{{\rm UT}}
\def\oCG{\overline{\CG}}
\def\hggu{\widehat{\ggu}}

\def\mJ{\mathbb{J}}

\begin{document}

\author{V. V. \  Bavula  
}

\title{Lie algebras of  triangular polynomial derivations and an isomorphism criterion for their Lie factor algebras}

\maketitle

\begin{abstract}
The Lie algebras $\ggu_n$ ($n\geq 2$) of triangular polynomial
derivations, their injective limit $\ggu_\infty$ and the completion $\hggu_\infty$  are studied in detail.
 The ideals of $\ggu_n$ are classified,  all
 of them are characteristic ideals. Using the classification of ideals, an  explicit  criterion is given for the Lie factor algebras of $\ggu_n$ and $\ggu_m$ to be isomorphic.  For (Lie)
algebras (and their modules) two new dimensions are introduced: the {\em central}
dimension $\cdim$ and the {\em uniserial} dimension $\udim$. 
 It is shown that $\cdim (\ggu_n) = \udim
(\ggu_n) = \o^{n-1}+\o^{n-2}+\cdots + \o +1$ for all $n\geq 2$ where $\o$ is the first infinite ordinal. Similar results are proved for the Lie algebras $\ggu_\infty$ and $\hggu_\infty$. In particular, $\udim (\ggu_\infty ) = \o^\o$ and $\cdim (\ggu_\infty ) =0$.

$\noindent $



$${\bf Contents}$$
\begin{enumerate}
\item Introduction.
 \item The derived and  upper  central series for the Lie algebra
$\ggu_n$.

\item Classification of ideals of the Lie algebra $\ggu_n$.

\item  The Lie algebras $\ggu_n$ are locally finite dimensional and locally nilpotent.

\item The  isomorphism problem for the factor algebras of the Lie  algebras $\ggu_n$.

\item The  Lie algebra $\ggu_\infty$.
\item The  Lie algebra $\hggu_\infty$.

\end{enumerate}
\end{abstract}


\section{Introduction}

Throughout,
module means
a left module;
 $\N :=\{0, 1, \ldots \}$ is the set of natural numbers; $K$ is a
field of characteristic zero and  $K^*$ is its group of units;
$P_n:= K[x_1, \ldots , x_n]=\bigoplus_{\alpha \in \N^n}
Kx^{\alpha}$ is a polynomial algebra over $K$ where
$x^{\alpha}:=x_1^{\alpha_1}\cdots x_n^{\alpha_n}$;
$\der_1:=\frac{\der}{\der x_1}, \ldots , \der_n:=\frac{\der}{\der
x_n}$ are the partial derivatives ($K$-linear derivations) of
$P_n$; $\Der_K(P_n) =\bigoplus_{i=1}^nP_n\der_i$ is the Lie
algebra of $K$-derivations of $P_n$; $A_n:= K \langle x_1, \ldots
, x_n , \der_1, \ldots , \der_n\rangle  =\bigoplus_{\alpha , \beta
\in \N^n} Kx^\alpha \der^\beta$  is  the $n$'th {\em Weyl
algebra}; for each natural number $n\geq 2$,
$$\ggu_n :=
K\der_1+P_1\der_2+\cdots +P_{n-1}\der_n$$ is the  {\em Lie algebra of triangular polynomial derivations} (it is a Lie subalgebra of the Lie algebra $\Der_K(P_n)$), $U_n:= U(\ggu_n)$ is
its universal enveloping algebra.


\begin{itemize}
\item (Proposition \ref{a8Dec11}) {\em The Lie algebras $\ggu_n$
where  $ n\geq 2$ are pairwise non-isomorphic. The Lie algebra
$\ggu_n$ is a solvable but not nilpotent Lie algebra. The centre
$Z(\ggu_n)$ of the Lie algebra $\ggu_n$ is $K\der_n$. All the
inner derivations of the Lie algebra $\ggu_n$ are locally
nilpotent derivations. The derived and upper central series are
found  for the Lie algebra $\ggu_n$.}
\end{itemize}

A new dimension for algebras and modules is introduced -- the {\em uniserial dimension} (Section \ref{CILAU}) -- that turned out to be  a very useful tool in studying non-Noetherian Lie algebras, their ideals and automorphisms, \cite{Lie-Un-AUT,Lie-Un-MON, Bav-Bod}.

\begin{itemize}

\item (Theorem \ref{10Dec11}) {\em  The Lie algebra
$\ggu_n$ is a uniserial, Artinian but not Noetherian Lie algebra
and its uniserial dimension is equal to $\udim (\ggu_n)  =
\o^{n-1}+\o^{n-2}+\cdots +\o +1$ where $\o$ is the first infinite
ordinal. A classification of all the ideals of $\ggu_n$  is given
and for each ideal an explicit basis is found. }
\item (Proposition \ref{b24Jan12}) {\em For every ideal of the Lie algebra $\ggu_n$, an explicit $K$-basis for its centralizer is found. The set of the centralizers of all the ideals of the Lie algebra $\ggu_n$ is a finite set that contains precisely 2n-1 elements.}
\end{itemize}

A Lie algebra $\CG$ is called a {\em locally nilpotent} ( respectively, {\em locally finite dimensional}) Lie algebra if every finitely generated Lie subalgebra of $\CG$ is a nilpotent (respectively, finite dimensional) Lie algebra.

\begin{itemize}
\item (Theorem \ref{22Dec11}) {\em The Lie algebras $\ggu_n$ are locally finite dimensional and locally nilpotent Lie algebras}.

 \item (Corollary  \ref{b10Dec11}) {\em  The ideal $\ggu_{n,2}:= P_1\der_2+
 \cdots + P_{n-1}\der_n$ is the largest proper
 ideal of the Lie algebra $\ggu_n$ and the centre $Z(\ggu_n) =
K\der_n$ of the Lie algebra $\ggu_n$ is the least nonzero one. The ideals $I_s:= \sum_{i=0}^{s-1}
Kx_1^i\der_n$ (where $s=1,2,\ldots $) are the only finite
dimensional ideals of the Lie algebra $\ggu_n$, and $\dim_K(I_s) =
s$. }
 \item (Corollary  \ref{c10Dec11}) {\em All the ideals of the Lie
 algebra $\ggu_n$ are characteristic ideals (i.e.,  invariant
 under the group of automorphisms of the Lie algebra $\ggu_n$). }
\item (Theorem \ref{A10Dec11}) {\em The central series $\{ Z^{(\l
)} (\ggu_n) \}_{\l \in \mW}$ stabilizers precisely on $\ord (\O_n)
:= \o^{n-1}+\o^{n-2}+\cdots + \o +1$ step, i.e.,  the central
dimension of the Lie algebra $\ggu_n$ is  $\cdim (\ggu_n) = \ord
(\O_n)$. Moreover, for each $\l \in [1,\ord (\O_n)]$, $Z^{(\l
)}(\ggu_n) = I_\l$ where $I_\l$ is given explicitly by
(\ref{Jlid}). }
\end{itemize}
In Section \ref{ISOPRUN}, Theorem \ref{24Dec11} and Corollary \ref{a24Dec11} give an answer to the question:

$\noindent $

{\em Let $I$ and $J$ be ideals of the Lie algebras $\ggu_n$ and $\ggu_m$ respectively. When are the Lie factor algebras $\ggu_n / I$ and $\ggu_m / J$  isomorphic?}

$\noindent $

 The answer is given in explicit terms (via the uniform dimensions of the ideals $I$ and $J$) using the classification of the ideals of the Lie algebras $\ggu_n$ (Theorem \ref{10Dec11}). In particular, there are only countably many ideals $I$ of the Lie algebra $\ggu_n$ such that $\ggu_n/I\simeq \ggu_n$. Theorem \ref{10Dec11}.(1)  shows that every ideal $I$ of the Lie algebra $\ggu_n$ is {\em uniquely determined} by its uniform dimension
 $$\l = \udim (I) \in [0, \o^{n-1}+\o^{n-2}+\cdots + \o +1],$$
  i.e., $I=I_\l$,  where $\o$ is the first infinite ordinal  number.

\begin{itemize}
\item {\rm (Corollary \ref{a24Dec11})}
{\em Let $n$ and $m$ be natural numbers such that $2\leq n <m$, $I$ and $J$ be ideals of  the Lie algebras $\ggu_n$ and $\ggu_m$ respectively. Then the Lie algebras $\ggu_n / I$ and $\ggu_m / J$ are isomorphic iff}
\begin{eqnarray*}
 (I,J)&\in &  \{ (I_\l , I_\mu ) \, | \, \l =i\o^{n-2}+\nu, \mu =\o^{m-1}+\o^{m-2}+\cdots +\o^n+j\o^{n-2}+\nu \\
 & &  where\;\; \nu\in [ 1, \o^{n-2})\cup \{ 0\} \;\;  and\;\; i,j\in \N\}\\
 &\cup & \bigcup_{s=2}^{n-1}\{ (I_\l , I_\mu ) \, | \, \l = \o^{n-1}+\o^{n-2}+\cdots +\o^s+ i\o^{s-2}+\nu, \mu =\o^{m-1}+\o^{m-2}+\cdots \\
  & &+\o^s +j\o^{s-2}+\nu\;\;   where\;\; \nu\in [ 1, \o^{s-2})\cup \{ 0\}\;\;  and\;\; i,j\in \N\}\\
  &\cup& \{ (I_\l , I_\mu )\, | \, \l = \o^{n-1}+\o^{n-2}+\cdots +\o+ \varepsilon , \mu =\o^{m-1}+\o^{m-2}+\cdots +\o +\varepsilon  \\
   & &   where\;\;  \varepsilon =0,1\}.
  \end{eqnarray*}
\end{itemize}

In Section \ref{LIUIN}, the Lie algebra $\ggu_\infty:= \cup_{n\geq 2}\ggu_n = \oplus_{n\geq 2} P_{n-1}\der_n$ is studied. Many of the properties of the Lie algebra $\ggu_\infty$ are similar to those of $\ggu_n$ ($n\geq 2$) but there are several differences. For example, the Lie algebra $\ggu_\infty$ is not solvable, not Artinian but almost Artinian, $\udim (\ggu_\infty ) = \o^\o$. A classification of all the ideals of the Lie algebra $\ggu_\infty$ is obtained (Theorem \ref{26Dec11}).  All  the ideals of the Lie algebra $\ggu_\infty$  are characteristic ideals (Corollary \ref{a26Dec11}). Corollary \ref{ca24Dec11} is an explicit criterion of when two Lie  factor algebras  of $\ggu_\infty$ are isomorphic.

In Section \ref{HHLIU}, the topological Lie algebra $\hggu_\infty$, which is the completion of the Lie algebra $\ggu_\infty$, is studied. Its properties diverge further from those of the Lie algebras $\ggu_n$ ($n\geq 2$) and $\ggu_\infty$. All the closed and all the open ideals of the topological Lie algebra $\hggu_\infty$ are classified (Theorem \ref{20Jan12}.(1)).
\begin{itemize}
\item {\rm (Corollary \ref{a21Jan12})}
\begin{enumerate}
\item {\em The topological Lie algebra $\hggu_\infty$ is an open uniserial, closed uniserial, open almost Artinian and closed almost Artinian Lie algebra which is neither open nor closed Artinian and is neither open nor closed  Noetherian.}
\item {\em The uniserial dimensions of the sets of open and of closed ideals of the topological Lie algebra $\hggu_\infty$ coincide and are equal  to} $\udim (\hggu_\infty ) = \o^\o$.
\item {\em All the open/closed ideals of the Lie algebra $\hggu_\infty$ are topologically characteristic ideals.}
\end{enumerate}
\end{itemize}


\section{The derived and  upper central series for the Lie algebra
$\ggu_n$}\label{TULSS}

In this section, various properties of the Lie algebras $\ggu_n$ are proved (Proposition \ref{a8Dec11}, Corollary \ref{b8Dec11}) that are used widely in the rest of the paper. At the end of the section, the image and the kernel  of the algebra homomorphism $\chi_n : U_n\ra A_{n-1}\t K[\der_n]$ (see (\ref{UnAn})) are found (Theorem \ref{12Jan12}). In particular, it is shown that the algebra $\chi_n (U_n)$ is not finitely generated and neither left nor right Noetherian.

Let $\CG$ be a Lie algebra over the field $K$ and $\ga$, $\gb$ be
its ideals. The {\em commutant}  $[\ga , \gb ]$ of the ideals
$\ga$ and $\gb$ is the linear span in $\CG$ of all the elements
$[a,b]$ where $a\in \ga$ and $b\in \gb$. The commutant $[U,V]$ of
  subspaces $U$ and $V$ of $\CG$ is defined in the same manner.
The commutant $[\ga , \gb ]$ is an ideal of the Lie algebra $\CG$
such that $[\ga , \gb ]\subseteq \ga \cap \gb$. In particular,
$\CG_{(1)}:= \CG^{(1)}:=[\CG , \CG ]$ is called the {\em
commutant} of the Lie algebra $\CG$. Let us define recursively the
following set of ideals of the Lie algebra $\CG$,
$$ \CG_{(i)}:=[\CG_{(i-1)}, \CG_{(i-1)}]\;\;\; {\rm and}\;\;\;
\CG^{(i)}:=[\CG , \CG^{(i-1)}]\;\;\; i\geq 2.$$ Clearly,
$\CG_{(i)}\subseteq \CG^{(i)}$ for all $i\geq 1$. The descending
chains of ideals of the Lie algebra $\CG$,
\begin{eqnarray*}
 \CG_{(0)}:=\CG &\supseteq &  \CG_{(1)}\supseteq \cdots \supseteq \CG_{(i)}\supseteq
   \CG_{(i+1)}\supseteq \cdots ,  \\
\CG^{(0)}:=\CG &\supseteq &  \CG^{(1)}\supseteq \cdots \supseteq
\CG^{(i)}\supseteq
   \CG^{(i+1)}\supseteq \cdots ,
\end{eqnarray*}
are called the {\em derived  series}  and the {\em upper central series} for the
Lie algebra $\CG$ respectively. Notice that
$$ \ggu_n=\bigoplus_{i=1}^n\bigoplus_{\alpha \in
\N^{i-1}}Kx^\alpha \der_i.$$ So, the elements 
\begin{equation}\label{Xai}
\Xai := x^\alpha \der_i = x_1^{\alpha_1}\cdots
x_{i-1}^{\alpha_{i-1}}\der_i, \;\;\; i=1, \ldots , n;  \; \alpha
\in \N^{i-1},
\end{equation}
form the $K$-basis $\CB_n$  for the Lie algebra $\ggu_n$. The basis
$\CB_n$  is called the {\em canonical basis} for $\ggu_n$.  For
all $1\leq i\leq j\leq n$, $\alpha \in \N^{i-1}$ and $\beta \in
\N^{j-1}$, 
\begin{equation}\label{Xai1}
[\Xai ,\Xbj]=\begin{cases}
0& \text{if }i=j,\\
\beta_iX_{\alpha +\beta -e_i, j}& \text{if }i<j,\\
\end{cases}
\end{equation}
where $e_1:=(1,0, \ldots , 0), \ldots , e_n:=(0,\ldots, 0, 1)$ is
the canonical free $\Z$-basis for the $\Z$-module $\Z^n$ and
$\N^i:= \sum_{k=1}^i\N e_k\subseteq \Z^i:= \sum_{k=1}^i\Z e_k$. In
particular, $\N \subseteq \N^2\subseteq \cdots \subseteq \N^n$ and
$\Z \subseteq \Z^2\subseteq \cdots \subseteq \Z^n$. The Lie
algebra $\ggu_n=\oplus_{i=1}^nP_{i-1}\der_i$ is the direct sum of
{\em abelian} (infinite dimensional when $i>1$)  Lie subalgebras $P_{i-1}\der_i$
(i.e.,  $[P_{i-1}\der_i, P_{i-1}\der_i]=0$) such that, for all $i<j$,
\begin{equation}\label{PdiPdj}
[P_{i-1}\der_i, P_{j-1}\der_j]=P_{j-1}\der_j.
\end{equation}
The inclusion ``$\subseteq $'' in (\ref{PdiPdj}) is obvious but
the equality follows from the fact that $[\der_i,
P_{j-1}]=P_{j-1}$. By (\ref{PdiPdj}), the Lie algebra $\ggu_n$
admits the finite, strictly descending chain of ideals
\begin{equation}\label{seruni}
\ggu_{n,1}:=\ggu_n\supset \ggu_{n,2}\supset \cdots \supset
\ggu_{n,i}\supset \cdots \supset \ggu_{n,n}\supset \ggu_{n,n+1}:=0
\end{equation}
where $\ggu_{n,i}:=\sum_{j=i}^nP_{j-1}\der_j$ for $i=1, \ldots ,
n$. By (\ref{PdiPdj}), for all $i<j$, 
\begin{equation}\label{cuni}
[\ggu_{n,i}, \ggu_{n,j}]\subseteq
\begin{cases}
\ggu_{n,i+1}& \text{if }i=j,\\
\ggu_{n,j}& \text{if }i<j.
\end{cases}
\end{equation}
For all $i=1, \ldots , n$, there is the canonical isomorphism of
Lie algebras 
\begin{equation}\label{cuni1}
\ggu_i\simeq \ggu_n/\ggu_{n,i+1}, \;\; X_{\alpha , j}\mapsto
X_{\alpha , j}+\ggu_{n, i+1}.
\end{equation}
 In particular, $\ggu_{n-1}\simeq \ggu_n/P_{n-1}\der_n$. Clearly,
$$ \ggu_2\subset \ggu_3\subset \cdots \subset \ggu_n\subset
\ggu_{n+1}\subset \cdots \subset \ggu_\infty :=\bigcup_{n\geq
2}\ggu_n=\bigoplus_{i\geq 1} \bigoplus_{\alpha \in \N^{i-1}}Kx^{\alpha}\der_i$$ is an ascending chain of Lie algebras. The polynomial
algebra $P_n$ is an $A_n$-module: for all elements $p\in P_n$,
$$ x_i*p=x_ip, \;\;\;\;\; \der_i*p = \frac{\der p}{\der x_i},
\;\;\; i=1, \ldots , n.$$ Clearly, $P_n\simeq A_n / \sum_{i=1}^n
A_n \der_i$, $1\mapsto 1+\sum_{i=1}^n A_n \der_i$. Since
$\ggu_n\subseteq A_n$, the polynomial algebra $P_n$ is also a
$\ggu_n$-module.

Let $V$ be a vector space over $K$. A $K$-linear map $\d : V\ra V$
is called a {\em locally nilpotent map} if $V=\cup_{i\geq 1} \ker
(\d^i)$ or, equivalently, for every $v\in V$, $\d^i (v) =0$ for
all $i\gg 1$. When  $\d$ is a locally nilpotent map in $V$ we
also say that $\d$ {\em acts locally nilpotently} on $V$. Every {\em nilpotent} linear map  $\d$, that is $\d^n=0$ for some $n\geq 1$, is a locally nilpotent map but not vice versa, in general.   Let
$\CG$ be a Lie algebra. Each element $a\in \CG$ determines the
derivation  of the Lie algebra $\CG$ by the rule $\ad (a) : \CG
\ra \CG$, $b\mapsto [a,b]$, which is called the {\em inner
derivation} associated with $a$. The set $\Inn (\CG )$ of all the
inner derivations of the Lie algebra $\CG$ is a Lie subalgebra of
the Lie algebra $(\End_K(\CG ), [\cdot , \cdot ])$ where $[f,g]:=
fg-gf$. There is the short exact sequence of Lie algebras
$$ 0\ra Z(\CG ) \ra \CG\stackrel{\ad}{\ra} \Inn (\CG )\ra 0,$$
that is $\Inn (\CG ) \simeq \CG / Z(\CG )$ where $Z(\CG )$ is the {\em centre} of the Lie algebra $\CG$ and $\ad ([a,b]) = [
\ad (a) , \ad (b)]$ for all elements $a, b \in \CG$. An element $a\in \CG$ is called a {\em locally nilpotent element} (respectively, a {\em nilpotent element}) if so is the inner derivation $\ad (a)$ of the Lie algebra $\CG$. Let $J$ be a non-empty subset of $\CG$ then
$\Cen_{\CG}(J) :=\{ a\in \CG \, | \, [a,b]=0$ for all $b\in J\}$
is called the {\em centralizer} of $J$ in $\CG$. It is a Lie subalgebra of the Lie algebra $\CG$. Let $A$ be an associative algebra and $I$ be a non-empty subset of
$A$. Then $\Cen_A(I):=\{ a\in A\, | \, ab=ba$ for all $b\in I\}$
is called the {\em centralizer} of $I$ in $A$. It is a subalgebra of $A$.

\begin{proposition}\label{a8Dec11}
\begin{enumerate}
\item The Lie algebra $\ggu_n$ is a solvable but not nilpotent Lie
algebra. \item The finite chain of ideals (\ref{seruni}) is the
 derived  series for the Lie algebra $\ggu_n$, that is
$(\ggu_n)_{(i)}= \ggu_{n,i+1}$ for all $i\geq 0$. \item The upper
central series for the Lie algebra $\ggu_n$ stabilizers at the
first step, that is  $(\ggu_n)^{(0)}=\ggu_n$ and
$(\ggu_n)^{(i)}=\ggu_{n, 2}$ for all $i\geq 1$. \item Each element
$u\in \ggu_n$ acts locally nilpotently on the $\ggu_n$-module
$P_n$. \item All the inner derivations of the Lie algebra $\ggu_n$
are locally nilpotent derivations.  \item The centre $Z(\ggu_n)$
of the Lie algebra $\ggu_n$ is $K\der_n$. \item The Lie algebras
$\ggu_n$ where  $ n\geq 2$ are pairwise non-isomorphic.
\end{enumerate}
\end{proposition}

{\it Proof}. 1. Statement 1 follows from statements 2 and 3.

2 and 3. Statements 2 and 3 follow from the decomposition
$\ggu_n=\oplus_{i=1}^n P_{i-1} \der_i$ and (\ref{PdiPdj}).

4. Statement 4 follows from the definition of the Lie algebra
$\ggu_n$.

5. Let $u\in \ggu_{n,i}\backslash \ggu_{n, i+1}$ for some $i=1,
\ldots , n$. Then $u=a\der_i+u'$ for some elements $a\in P_{i-1}$
and $u'\in \ggu_{n,i+1}$. Let $\d=\ad (u)$,   $\der = \ad
(a\der_i)$, and  $v\in \ggu_n$. We have to show that $\d^s (v) =0$
for all $s\gg 1$. By applying (\ref{cuni}) twice,  we see that $\d
(v) \in \ggu_{n,i}$ and $ \d^2 (v) \in \ggu_{n,i+1}$. By replacing
the element $v$ with $\d^2 (v)$, without loss of generality we may
assume that $v\in \ggu_{n,i+1}$ and $v\neq 0$. In view of the
equality $\ggu_{n, i+1} = \bigoplus_{j=i+1}^n P_{j-1}\der_j$,
$v=b\der_{i+1}+v'$ for some elements  $0\neq b\in P_i$ and $v'\in
\ggu_{n, i+2}$. For all natural numbers $s\geq 1$, by
(\ref{cuni}),
$$ \d^s (v) \equiv \d^s(b\der_{i+1}) \equiv \d^s (b) \der_{i+1} \equiv
\der^s(b)\der_{i+1} \mod \ggu_{n,i+2}.$$ By statement 4, $\der^s(b)=0$
for all $s\gg 1$. Then $\d^s (v) \in \ggu_{n,i+2}$ for all $s\gg
1$. Similarly, $\d^s(v) \in \ggu_{n,i+3}$ for all $s\gg 1$.
Applying the same argument several times we see that
$ \d^s (v)\in \ggu_{n,n+1}=0$   for all $s\gg
1$.  This means that $\d$ is a locally nilpotent map, as required.

6. It is well-known and easy to show that 
\begin{equation}\label{CenAnd}
\Cen_{A_n}(\der_1, \ldots , \der_n )= K[\der_1, \ldots , \der_n].
\end{equation}
It follows that  
\begin{equation}\label{CenAnd1}
\Cen_{A_n}(\der_1, \ldots , \der_n, x_1\der_n, \ldots ,
x_{n-1}\der_n )= K[\der_n],
\end{equation}
and then
$$Z(\ggu_n) \subseteq \ggu_n\cap \Cen_{A_n}(\der_1, \ldots , \der_n, x_1\der_n, \ldots ,
x_{n-1}\der_n )=\ggu_n\cap  K[\der_n]=K\der_n.$$ The reverse
inclusion, $Z(\ggu_n) \supseteq K\der_n$, is obvious. Therefore,
$Z(\ggu_n) = K\der_n$.

7. Statement 7 follows from statement 2 as the derived
series for the algebras $\ggu_n$ have distinct lengths (and they
are isomorphism invariants).  $\Box $


Proposition \ref{a8Dec11}.(5) allows us to produces many automorphisms of the Lie algebra $\ggu_n$. For every element $a\in \ggu_n$, the inner derivation $\ad (a)$ is a locally nilpotent derivation, hence $$e^{\ad (a)}:=\sum_{i\geq 0} \frac{\ad (a)^i}{i!}\in \Aut_K(\ggu_n).$$ In \cite{Lie-Un-AUT}, the group $\Aut_K(\ggu_n)$ of automorphisms of the Lie algebra $\ggu_n$ and its explicit generators are found and it was shown that the {\em adjoint  group} $\langle e^{\ad (a)}\, | \, a\in \ggu_n\rangle$ is a {\em tiny} part of the group $\Aut_K(\ggu_n)$.

The next lemma classifies the {\em nilpotent} inner derivations of the Lie algebras $\ggu_n$.

\begin{lemma}\label{b19Jan12}
Let $a\in \ggu_n$ and $\d =\ad (a)$. The following statements are equivalent.
\begin{enumerate}
\item The map $\d$  is a nilpotent derivation of the Lie algebra $\ggu_n$.
\item $a\in P_{n-1}\der_n$.
\item $\d^2=0$.
\end{enumerate}
\end{lemma}

{\it Proof}.  The implications $(2\Rightarrow 3\Rightarrow 1)$ are obvious.

 $(1\Rightarrow 2)$ We have to show that the derivation $\d$ is not nilpotent for every element $a\in \ggu_n\backslash P_{n-1}\der_n$. Let $u\in \ggu_n\backslash P_{n-1}\der_n$, i.e.,  $u=p_i\der_i+p_{i+1}\der_{i+1}+\cdots +p_n\der_n= p_i\der_i+\cdots$ where $p_j\in P_{j-1}$ for all $j=i, \ldots , n$, $i<n$ and $p_i\neq 0$. Since
 $$ \d^m (x_i^m\der_{i+1})=m!p_i^m\der_{i+1}+\cdots\;\;\; {\rm for \; all}\;\; m\geq 1,$$the derivation $\d$ is not a nilpotent derivation.  $\Box $


Let $A$ be a ring. A subset $S$ of $A$ is called a {\em
multiplicative subset} or a {\em multiplicatively closed subset}
of $A$ if $1\in S$, $SS\subseteq S$ and $0\not\in S$. Every
associative algebra $A$ can be seen as a Lie algebra $(A, [\cdot,
\cdot ])$ where $[a,b]=ab-ba$ is the {\em commutator} of elements
$a,b\in A$. For each element $a\in A$, the map $\ad (a): A\ra A$,
$b\mapsto [a,b]$, is a $K$-derivation of the algebra $A$ seen  as
an associative and  Lie algebra. The derivation $\ad (a)$ is
called the {\em inner derivation} of $A$ associated with the
element $a$. So, the associative algebra $A$ and the Lie algebra
$(A, [\cdot , \cdot ])$ have the same set of inner derivations
$\Inn (A)$ and the same centre $Z(A)$.

Let $\d$ be a derivation of a ring $A$. For all elements $a,b\in A$,
\begin{equation}\label{dnab}
\d^n(ab) = \sum_{i=0}^n {n\choose i} \d^i(a)\d^{n-i}(b), \;\;\;\; n\geq 1;
\end{equation}
\begin{equation}\label{dnab1}
a^nb = \sum_{i=0}^n {n\choose i} (\ad\, a)^i(b)a^{n-i}, \;\;\;\; n\geq 1.
\end{equation}

\begin{corollary}\label{b8Dec11}

\begin{enumerate}
\item  The inner derivations $\{ \ad (u)\, | \, u\in \ggu_n\}$ of the universal enveloping
algebra $U_n$ of the Lie algebra $\ggu_n$  are locally nilpotent
derivations. \item Every multiplicative subset $S$ of $U_n$ which is generated by an arbitrary set of elements of $\ggu_n$ is a
(left and right) Ore set in $U_n$. Therefore, $S^{-1}U_n\simeq
U_nS^{-1}$.
\end{enumerate}
\end{corollary}

{\it Proof}.  1. Statements 1 follows from (\ref{dnab}) and
Proposition \ref{a8Dec11}.(5).

2. Statements  2 is an  easy corollary of (\ref{dnab1}) and
Proposition \ref{a8Dec11}.(5).  $\Box $


{\it Example}. The set $S=\{ \der^\alpha \, | \, \alpha\in \N^n\}$
is a multiplicative subset of the algebra $U_n$. By Corollary
\ref{b8Dec11}.(2), the localization ring $S^{-1}U_n$ exists.

{\bf Generalized Weyl Algebras}. Let $D$ be a ring, $\sigma=(\sigma_1,...,\sigma_n)$
 an $n$-tuple  of  commuting automorphisms of $D$,
($\sigma_i\sigma_j=\sigma_j\sigma_i$, for
all $i,j$), and $a=(a_1,...,a_n)$  an $n$-tuple  of (non-zero) elements of the centre $Z(D)$  of $D$, such that $\sigma_i(a_j)=a_j$ for all $i\neq j$.

        The {\it generalized Weyl algebra} $A=D(\sigma,a)$ (briefly GWA) of
degree $n$ with {\it base } ring  $D$  is  the  ring  generated
by $D$  and    $2n$ indetermina\-tes $x_1,...,x_n,$ $y_1,...,$ $y_n$
subject to the defining relations \cite{Bav-GWA-FA-91}, \cite{Bav-GWA-AA}:
$$y_ix_i=a_i,\qquad x_iy_i=\sigma_i(a_i),$$
$$x_i\alpha=\sigma_i(\alpha)x_i,\;\; {\rm and }\;\;
 y_i\alpha=\sigma_i^{-1}(\alpha)y_i,\;\; \;\; {\rm for }\;\;   \alpha \in D,$$
$$[x_i,x_j]=[y_i,y_j]=[x_i,y_j]=0, \;\;{\rm for }\;\;  i\neq j,$$
where $[x, y]=xy-yx$. We say that  $a$  and $\sigma $ are the
sets  of {\it defining } elements and automorphisms of $A$
respectively. For a vector $k=(k_1,...,k_n)\in \mathbb{Z}^n$ we
put $v_k=v_{k_1}(1)\cdots v_{k_n}(n)$ where, for  $1\leq i\leq n$
and $m\geq 0$:  $\,v_{m}(i)=x_i^m$, $\,v_{-m}(i)=y_i^m$,
$v_0(i)=1$. It follows from  the  definition  of the  GWA that
$$A=\oplus_{k\in \mathbb{Z}^n} A_k$$
is a $\mathbb{Z}^n$-graded algebra ($A_kA_e\subset A_{k+e},$ for
all  $k,e \in \mathbb{Z}^n$),
 where $A_k=Dv_k=v_kD$.

Let $\CP_n$ be the polynomial algebra  $K[H_1, \ldots , H_n]$ in
$n$ indeterminates  and let $\sigma=(\sigma_1,...,\sigma_n)$ be an
$n$-tuple  of  commuting automorphisms of $\CP_n$ defined as
follows: $\s_i(H_i)=H_i-1$ and $\s_i(H_j)=H_j$, for $i\neq j$. The
map
$$ A_n\ra \CP_n ((\sigma_1,...,\sigma_n), (H_1, \ldots , H_n)), \;\;
x_i\mapsto x_i, \; \; \der_i\mapsto  y_i, \; \; {\rm for }\;i=1, \ldots ,
n, $$ is a $K$-algebra isomorphism. We identify the Weyl algebra
$A_n$ with the GWA above via this isomorphism.
The Weyl algebra $A_n=\oplus_{\alpha \in \Z^n} A_{n, \alpha}$ is a $\Z$-graded algebra ($A_{n,\alpha}A_{n, \beta}\subseteq A_{n,\alpha +\beta }$ for all $\alpha , \beta \in \Z^n$).

  The multiplicative sets $S=\{ \der^\alpha \, | \, \alpha\in \N^n\}$ and $T=\{ x^\alpha \, | \, \alpha\in \N^n\}$ are (left and right) Ore sets of the Weyl algebra $A_n$ and
\begin{eqnarray*}
B_n&:=&S^{-1}A_n =K[H_1, \ldots , H_n][\der_1^{\pm 1}, \ldots , \der_n^{\pm 1}; \s_1^{-1}, \ldots , \s_n^{-1}],  \\
B_n'&:=& T^{-1}A_n = K[H_1, \ldots , H_n][x_1^{\pm 1},
\ldots , x_n^{\pm 1}; \s_1, \ldots , \s_n ].
\end{eqnarray*}
 The Weyl algebra $(A_n,[\cdot ,
\cdot ])$ is a $\Z^n$-graded Lie algebra, that is $[A_{n,\alpha},
A_{n,\beta}]\subseteq A_{n, \alpha +\beta}$ for all elements
$\alpha, \beta \in \Z^n$. By the very definition, $\ggu_n$ and
$\Der_K(P_n)$ are $\Z^n$-graded Lie subalgebras of the Lie algebra
$A_n$.

The Lie algebra $\ggu_n$ contains both finite and infinite  dimensional maximal abelian Lie subalgebras as the next lemma shows.
\begin{lemma}\label{e10Dec11}
\begin{enumerate}
\item  $\CD_n:=\bigoplus_{i=1}^n K\der_i$ is a maximal abelian Lie
subalgebra of the Lie algebra $\ggu_n$ and $\Cen_{\ggu_n}(\CD_n) = \CD_n$. \item The ideal
$P_{n-1}\der_n$ of the Lie algebra $\ggu_n$ is a maximal abelian
Lie subalgebra of $\ggu_n$ and $\Cen_{\ggu_n}(P_{n-1}\der_n) = P_{n-1}\der_n$.
\end{enumerate}
\end{lemma}

{\it Proof}. 1. Statement 1 follows from the following two facts
$\Cen_{A_n}(\CD_n) = K[\der_1, \ldots , \der_n]$ (see (\ref{CenAnd})) and $\Cen_{\ggu_n}(\CD_n)=\ggu_n\cap
K[\der_1, \ldots , \der_n]=\CD_n$.

2. It follows from the equality $\Cen_{A_{n-1}}(P_{n-1})=P_{n-1}$
that
$$\Cen_{A_{n-1}\t K[\der_n]}(P_{n-1}\der_n)
=\Cen_{A_{n-1}}(P_{n-1}) \t K[\der_n]=P_{n-1}[\der_n].$$ Then
$\Cen_{\ggu_n}(P_{n-1}\der_n) = \ggu_n\cap P_{n-1} [ \der_n]=
P_{n-1} \der_n$. Therefore, the ideal $P_{n-1}\der_n$ is a maximal
abelian  Lie subalgebra  of $\ggu_n$. $\Box $


{\bf The homomorphism $\chi_n$}. The inclusion $\ggu_n\subseteq A_{n-1}\t K[\der_n]$, induces the algebra
homomorphism 
\begin{equation}\label{UnAn}
\chi_n : U_n\ra A_{n-1}\t K[\der_n], \;\; \Xai \mapsto x^\alpha \der_i.
\end{equation}
The image and the kernel of $\chi_n$ are found in Theorem \ref{12Jan12}. By Corollary \ref{b8Dec11}.(2),  the homomorphism $\chi_n$ can be extended to the algebra
homomorphism (where $S=\{ \der^\alpha \, | \, \alpha \in \N^\alpha \}$)
$$ \chi_n : S^{-1}U_n\ra S^{-1}(A_{n-1}\t K[\der_n])=B_{n-1}\t K[\der_n , \der_n^{-1}].$$
It is obvious that $\chi_n (S^{-1}U_n) = B_{n-1}\t K[\der_n ,
\der_n^{-1}]$ (since $X_{\alpha , i}\der_i^{-1}\mapsto x^\alpha \der_i\der_i^{-1}=x^\alpha $ and $x_i=\der_i^{-1}H_i$).

Let us define the relation $\prec$ on the set $\N^n$: we write $\alpha\prec \beta$ for elements $\alpha = (\alpha_i)$ and $ \beta = (\beta_i)$ of $\N^n$ iff either $\alpha =0$ and $\beta$ is arbitrary (i.e.,  $0\prec\beta$ for all $\beta \in \N^n$) or $\alpha\neq 0$, $\beta \neq 0$ and $\max \{ i \, | \, \alpha_i\neq 0\}<\max \{ i \, | \, \beta_i\neq 0\}$.  Clearly, $\alpha\prec \beta$ and $\beta \prec \g$ imply $\alpha\prec \g$; $\alpha \prec \alpha$ iff $\alpha =0$; for all $\alpha , \beta \in \N^n\backslash \{ 0\}$ , $\alpha\prec \beta$ implies $\beta\not\prec\alpha$.

\begin{theorem}\label{12Jan12}

\begin{enumerate}
\item The set $W_n:=\{ x^\alpha\der^\beta\, | \, \alpha, \beta \in \N^n; \alpha\prec \beta\}$ is a $K$-basis for the algebra $U_n':=\chi_n (U_n)$.
\item The kernel of the algebra homomorphism $\chi_n $ is the ideal of the algebra $U_n$ generated by the elements
    $X_{\alpha , i}X_{\beta , j}- X_{0, i} X_{\alpha + \beta, j}$ where $i\leq j$,  $\alpha\in \N^{i-1}$ and $\beta\in \N^{j-1}$.
    \item The homomorphism $\chi_n$ is not a surjective map.

        \item The algebra $U_n'$ is not a finitely generated, neither left nor right Noetherian algebra.

\end{enumerate}
\end{theorem}

{\it Proof}. 1. The elements of  the set $W_n$  are
 $K$-linearly independent elements since they are so as elements  of the algebra $A_{n-1}\t K[\der_n]$. Let $\CI_n:=\sum_{w\in W_n}Kw$. We have to show that $\chi_n(U_n) = \CI_n$. The algebra $\chi_n(U_n)$  is generated  by the elements $\chi_n(X_{\alpha , i}) = x^\alpha \der_i$.  Using the relations
 $$ x^\alpha\der_ix^\beta \der_i = x^{\alpha +\beta } \der_i^2, \;\; [x^\alpha\der_i,x^{\g} \der_j] = \g_i x^{\alpha +\g - e_i} \der_i\der_j\;\; {\rm for } \;\; i<j, $$ we see that the algebra $\chi_n(U_n)$ is contained in the linear span, say $\CI_n'$, of the elements
 $$ \der_1^{\beta_1}x^{\alpha_2}\der_2^{\beta_2}\cdots x^{\alpha_i}\der_i^{\beta_i}\cdots x^{\alpha_n}\der_n^{\beta_n}\;\; {\rm where} \;\; \beta_i\in \N, \; \alpha_i\in \N^{i-1}, \; i=2, \ldots , n.$$ Using the commutation relations $[\der_i, x_j]=\d_{ij}$ (where $\d_{ij}$ is the Kronecker delta-function), every such element can be written as a linear combination of some elements of $W_n$. Therefore, $\chi_n(U_n) \subseteq \CI_n'\subseteq \CI_n$.

 To prove that the reverse inclusion $\chi_n(U_n)\supseteq \CI_n$ holds and hence to finish the proof of statement 1, we have to show that each element $w\neq 1$ of $W_n$ belongs to the algebra $\chi_n(U_n)$. The element $w\neq 1$ is the product
 $$ x_1^{\alpha_1}\cdots x_s^{\alpha_s}\der_1^{\beta_1} \cdots \der_t^{\beta_t}\;\; {\rm for \; some}\;\; \alpha_i, \beta_j\in \N, s<t\; {\rm and}\; \beta_t\neq 0.$$
 The case where $\alpha_1=\cdots = \alpha_s=0$ is obvious. So, we can assume that $\alpha_s\neq 0$.  Each element $x_1^{\alpha_1}\cdots x_s^{\alpha_s}\der_1^{\beta_1}\cdots \der_s^{\beta_s}$ of the Weyl algebra $A_s$ can be written as a sum $\sum_{u,v\in \N^s} \l_{uv}\der^ux^v$ where $\l_{uv}\in K$. Now,
 $$w=(\sum_{u,v\in \N^s} \l_{uv}\der^ux^v)\der_{s+1}^{\beta_{s+1}}\cdots \der_t^{\beta_t}=
 \sum_{u,v\in \N^s} \l_{uv}\der^u \der_{s+1}^{\beta_{s+1}}\cdots \der_{t-1}^{\beta_{t-1}}  (x^v\der_t^{\beta_t})\in \chi_n(U_n).$$ The proof of statement 1 is complete. Moreover, the last step implies that the set
\begin{equation}\label{Wpn}
W_n':=\{ \der^{\alpha}, \der^{\beta} x^\nu \der_t^i\, | \, \alpha \in \N^n; i\geq 1; \beta , \nu \in \N^{t-1}\; \; {\rm and}\;\; \nu \neq 0;\;  t=2, \ldots , n\}
\end{equation}
 is  also a $K$-basis for the algebra $\chi_n (U_n)$. In more detail, the algebra $\chi_n(U_n)$ is a linear span of $W_n'$ and the elements of the set  $W_n'$ are linear independent since they are so as elements of the Weyl algebra $A_n$. Therefore, $W_n'$ is a $K$-basis of the algebra $\chi_n(U_n)$.  This basis is used in the proof of statement 2.

2. Let $I$ be the ideal of the algebra $U_n$ generated by the elements (the would be generators for $\ker (\chi_n)$) of statement 2. For each element $w'$ of the set $W_n'$, that is for $\der^\alpha$ and $\der^\beta x^v\der_t^i$, choose the element $w''\in \chi_n^{-1}(w')$ as follows
$$ w'':= \begin{cases}
\prod X_{0,i}^{\alpha_i}& \text{if }w'=\der^\alpha, \alpha = (\alpha_i),\\
\prod X_{0,i}^{\beta_i}\cdot X_{v,t}\cdot X_{0,t}^{i-1}& \text{if }w'=\der^\beta x^v\der_t^i.\\
\end{cases}
$$
So, $\chi_n(w'') = w'$ for all elements $w'\in W_n'$. Let $W_n''$ be the set of all elements $w''$. The elements in $W_n''$ are linearly independent in the algebra $U_n$ as pre-images of linearly independent elements. Let $\CI_n'':= \sum_{w''\in W_n''}Kw''$. To finish the proof of statement 2, it suffices to show that
$$U_n=\CI_n''+I.$$ (Suppose that the equality holds. Since $I\subseteq \ker (\chi_n)$ and the set $W_n''$ is mapped bijectively onto the basis $W_n'$ of the algebra $\chi_n(U_n)$, these two facts necessarily imply that $I=\ker (\chi_n)$). To prove the equality $U_n = \CI_n''+I$ we follow the line of the proof of statement 1. The relations (\ref{Xai1}) and $X_{\alpha , i}X_{\beta , i}\equiv X_{0,i} X_{\alpha +\beta , i} \mod I$ imply that the linear span of the elements
$$ \{ X_{0,1}^{\beta_1}(X_{0,2}^{\beta_2}X_{\alpha_2,2}')\cdots  (X_{0,i}^{\beta_i}X_{\alpha_i,i}') \cdots (X_{0,n}^{\beta_n}X_{\alpha_n,n}')\, | \, \beta_i\in \N, \alpha_i\in \N^{i-1}, X_{\alpha_i, i}'\in \{ 1, X_{\alpha_i, i}\}, \; i=2, \ldots , n\},$$
where
$$(X_{0,i}^{\beta_i}X_{\alpha_i,i}')=\begin{cases}
1& \text{if }\beta_i=0,\\
X_{0,i}^{\beta_i}X_{\alpha_i,i}'& \text{if }\beta_i\neq 0,\\
\end{cases}$$ generate the algebra $U_n$ modulo the ideal $I$. Using the generators for the  ideal $I$ and the relations (\ref{Xai1}), each of these elements can be written as a linear combination of the elements $w''$ (defined above). Then $U_n = \CI_n''+I$.

3. Suppose that the homomorphism $\chi_n$ is a surjective map, we seek a contradiction. The ideal $\ga := U_n\ggu_n U_n$ of the algebra $U_n$ contains $\ker (\chi_n)$, $U_n/\ga  \simeq K$, $\chi_n(\ga ) $ is an ideal of the algebra $\chi_n(U_n)=A_{n-1}\t K[\der_n]$ and
$$ A_{n-1}\t K[\der_n]/\chi_n (\ga ) = \chi_n(U_n) / \chi_n(\ga )= \chi_n(K+\ga ) / \chi_n(\ga ) = (K+\chi_n(\ga )) / \chi_n(\ga ) \simeq K.$$ The Weyl algebra $A_{n-1}$ is a simple infinite dimensional algebra, so it is mapped isomorphically onto its image under the algebra homomorphism
$$A_{n-1}\ra A_{n-1}\t K[\der_n]/\chi_n(\ga )\simeq K, \;\; a\mapsto a\t 1+\chi_n(\ga ),$$
 a contradiction.

4. See Proposition \ref{a14Jan12} where  a stronger statement is proved.  $\Box $


\begin{corollary}\label{x12Jan12}
The set $W_n'$ (see (\ref{Wpn})) is a $K$-basis for the algebra $U_n'$.
\end{corollary}

{\it Proof}. This fact was established in the proof of statement 1 of Theorem \ref{12Jan12}. $\Box $



\section{Classification of ideals of the Lie algebra $\ggu_n$}\label{CILAU}

In this section, the uniserial and central dimensions are
introduced and a classification of ideals is given for the Lie
algebra $\ggu_n$ (Theorem \ref{10Dec11}.(1)). It is proved that
the Lie algebra $\ggu_n$ is a uniserial,  Artinian but not
Noetherian  Lie algebra of uniserial dimension $\udim (\ggu_n)
=\o^{n-1}+\o^{n-2}+\cdots + \o +1$ (Theorem \ref{10Dec11}.(2));
every ideal of $\ggu_n$ is a characteristic ideal  (Corollary
\ref{c10Dec11}); the central series of the Lie algebra $\ggu_n$ is
found and it is shown that the central dimension of the Lie
algebra $\ggu_n$ is equal to $\cdim (\ggu_n)
=\o^{n-1}+\o^{n-2}+\cdots + \o +1$ (Theorem \ref{A10Dec11}).

{\bf The uniserial dimension}.  Let $(S, \leq )$ be a {partially
ordered set} (a {\em poset}, for short),  i.e.,   a set $S$ admits a
relation $\leq$ that satisfies three conditions: for all $a,b,c\in
S$,

(i) $a\leq a$;

(ii) $a\leq b$ and $b\leq a$ imply $a=b$;

(iii) $a\leq b$ and  $b\leq c$ imply $a\leq c$.

A poset $(S, \leq )$ is called an {\em Artinian} poset is every
non-empty subset $T$ of $S$ has a {\em minimal element}, say $t\in
T$, that is $t\leq t'$ for all $t'\in T$.  A poset $(S,\leq )$ is
a {\em well-ordered} if for all elements $a,b\in S$ either $a\leq
b$ or $b\leq a$. A bijection $f: S\ra S'$ between two posets $(S,
\leq )$ and $(S', \leq )$ is an {\em isomorphism} if $a\leq b$ in
$S$ implies $f(a)\leq f(b)$ in $S'$. Recall that the {\em ordinal
numbers} are the isomorphism classes of well-ordered Artinian
sets. The ordinal number (the isomorphism class) of a well-ordered
Artinian set $(S, \leq )$ is denoted by $\ord (S)$. The class of
all ordinal numbers is denoted by $ \mW $. The class $\mW$ is
well-ordered  by `inclusion' $\leq$ and Artinian. An associative
addition `$+$' and an associative multiplication `$\cdot$' are
defined in $\mW$ that extend the addition and multiplication of
the natural numbers.  Every non-zero natural number $n$ is
identified with $\ord (1<2<\cdots <n)$. Let $\o := \ord (\N , \leq
)$. More details on the ordinal numbers the reader can find in the book \cite{Rotman-Kneebone-Book}.

$\noindent $

{\it Definition}. Let $(S, \leq )$ be a partially ordered set. The {\em uniserial dimension} $\udim (S)$ of $S$ is the supremum of  $\ord (\CI )$ where $\CI$ runs through
all the Artinian well-ordered subsets of $S$.

$\noindent $

For a Lie algebra $\CG$, let $\CJ_0 (\CG)$ and $\CJ (\CG )$ be the
sets of all and all non-zero ideals  of the Lie algebra $\CG$,
respectively. So, $\CJ_0(\CG ) = \CJ (\CG )\cup \{ 0\}$. The sets
$\CJ_0(\CG)$ and $\CJ (\CG )$ are posets with respect to
inclusion. A Lie algebra $\CG$ is called {\em Artinian}
(respectively, {\em Noetherian}) if the poset $\CJ (\CG )$ is
Artinian (respectively, Noetherian). This means that every
descending (respectively, ascending) chains of ideals stabilizers.
A Lie algebra $\CG$ is called a {\em uniserial} Lie algebra if the
poset $\CJ (\CG )$ is a well-ordered set. This means that for any
two ideals  $\ga$ and $\gb$ of the Lie algebra $\CG$ either $\ga
\subseteq \gb$ or $\gb \subseteq \ga$.

$\noindent $

{\it Definition}. Let $\CG$ be an Artinian uniserial Lie algebra.
The ordinal number $\udim (\CG) := \ord (\CJ (\CG ))$ of the
Artinian  well-ordered set $\CJ (\CG )$ of nonzero ideals of $\CG$ is
called the {\em uniserial dimension} of the Lie algebra $\CG$. For
an arbitrary Lie algebra $\CG$, the uniserial dimension $\udim
(\CG )$ is the supremum of $\ord (\CI )$ where $\CI$ runs through
all the Artinian well-ordered sets of ideals.

$\noindent $

If $\CG$ is a Noetherian Lie algebra then $\udim (\CG ) \leq \o$.
So, the uniserial dimension is a measure of deviation from the
Noetherian condition. The concept of the uniserial dimension makes
sense for any class of algebras (associative, Jordan, etc.).

Let $A$ be an algebra and $M$ be its module, and let $\CJ_l(A)$ and $\CM (M)$ be the sets of all the nonzero left ideals of $A$ and of all the nonzero submodules of $M$, respectively. They are posets with respect to $\subseteq $. The {\em left uniserial dimension} of the algebra $A$ is defined as $\udim (A):= \udim (\CJ_l(A))$ and the {\em uniserial dimension} of the $A$-module $M$ is defined as $\udim (M):=\udim (\CM (M))$.

{\bf An Artinian well-ordering on the canonical basis $\CB_n$ of
$\ggu_n$}.  Let us define an Artinian well-ordering $\leq $ on the
canonical basis $\CB_n$  for the Lie algebra $\ggu_n$ by the rule:
$\Xai
>\Xbj$ iff $i<j$ or $i=j$ and $\alpha_{n-1}=\beta_{n-1}, \ldots,
\alpha_{m+1}=\beta_{m+1}, \alpha_m>\beta_m$ for some $m$.


{\it Examples}. For $n=2$, $ \der_2<x_1\der_2<x_1^2\der_2<\cdots
<\der_1 $.

 For $n=3$,
\begin{eqnarray*}
\der_3&< &x_1\der_3<x_1^2\der_3<\cdots<\\
x_2\der_3&< &x_1x_2\der_3<x_1^2x_2\der_3<\cdots < \\
&\cdots & \\
x_2^i\der_3&< &x_1x_2^i\der_3<x_1^2x_2^i\der_3<\cdots < \\
&\cdots & \\
\der_2&<&x_1\der_2<x_1^2\der_2<\cdots < \der_1.
\end{eqnarray*}
The next lemma is a straightforward consequence of the definition of the ordering $<$, we write $0<X_{\alpha , i}$ for all $X_{\alpha , i}\in \CB_n$.
\begin{lemma}\label{a18Dec11}
If $X_{\alpha , i}>X_{\beta , j}$ then
\begin{enumerate}
\item $X_{\alpha +\g, i}>X_{\beta +\g, j}$ for all $\g \in \N^{i-1}$,
\item $X_{\alpha -\g, i}>X_{\beta -\g, j}$ for all $\g \in \N^{i-1}$ such that $\alpha -\g \in \N^{i-1}$ and $\beta -\g \in \N^{j-1}$,
\item $[\der_k, X_{\alpha , i}]>[\der_k, X_{\beta , j}]$ for all $k=1, \ldots , i-1$ such that $\alpha_k\neq 0$, and
    \item $[X_{\g ,k}, X_{\alpha , i}]>[X_{\g ,k}, X_{\beta , j}]$ for all $X_{\g , k}>X_{\alpha , i}$ such that
           $[X_{\g ,k}, X_{\alpha , i}]\neq 0$, i.e.,  $\alpha_k\neq 0$.
\end{enumerate}
\end{lemma}

Let $\O_n$ be the set of indices $\{ (\alpha , i)\}$ that
parameterizes the canonical basis $\{ \Xai \}$ of the Lie algebra
$\ggu_n$. The set $(\O_n, \leq )$ is an Artinian well-ordered set,
where $(\alpha , i)\geq (\beta , j)$ iff $\Xai \geq \Xbj$, which
is isomorphic to the Artinian well-ordered set $(\CB_n, \leq )$
via $(\alpha , i)\mapsto \Xai$. We identify the posets $(\O_n,
\leq )$ and $(\CB_n , \leq )$ via this isomorphism.  It is obvious
that 
\begin{equation}\label{ordn}
\ord (\CB_n)=\ord (\O_n ) = \o^{n-1}+\o^{n-2}+\cdots +\o +1,
\end{equation}
$\O_2\subset \O_3\subset \cdots $ and $\CB_2\subset \CB_3\subset
\cdots $.   Let $[1, \ord (\O_n)]:=\{ \l \in \mW \, | \, 1\leq \l
\leq \ord (\O_n)\}$. By (\ref{Xai1}), if $[\Xai, \Xbj ]\neq 0$
then 
\begin{equation}\label{XaXj}
[\Xai , \Xbj ] <\min \{ \Xai , \Xbj \}.
\end{equation}
By (\ref{XaXj}), the  map 
\begin{equation}\label{Jlid}
\rho_n: [1, \ord (\O_n)]\ra \CJ (\ggu_n), \;\; \l  \mapsto I_\l :=I_{\l, n} :=
\bigoplus_{ (\alpha , i)\leq \l } K\Xai ,
\end{equation}
is a monomorphism  of posets ($\rho_n$  is an order-preserving
injection). We will prove that the map $\rho_n$ is a bijection
(Theorem \ref{10Dec11}.(1)) and as a result we will  have a
classification of all the ideals of the Lie algebra $\ggu_n$. Each
non-zero element $u$ of $\ggu_n$ is a finite linear combination
$$ u = \l\Xai +\mu \Xbj +\cdots +\nu X_{\s , k}= \l \Xai +\cdots $$
where $\l , \mu , \ldots , \nu \in K^*$ and $\Xai > \Xbj >\cdots
>X_{\s , k}$. The elements $\l \Xai$ and $\l \in K^*$ are called
the {\em leading term} and the {\em leading coefficient} of $u$
respectively, and the ordinal number denoted by $\ord (\Xai ) = \ord (\alpha
, i) \in [ 1, \ord (\O_n)]$, which is, by definition, the ordinal number that represents the Artinian well ordered set $\{ (\beta, j) \in \O_n \, | \, (\beta , j)\leq (\alpha , i)\}$,  is called the {\em ordinal degree} of
$u$ denoted by $\ord (u)$ (we hope that this notation will not
lead to confusion). For all non-zero elements $u,v\in \ggu_n$ and
$\l \in K^*$,

(i) $\ord (u+v) \leq \max \{ \ord (u) , \ord (v)\}$ provided
$u+v\neq 0$;

(ii) $\ord (\l u) = \ord (u)$;

(iii) $\ord ([u,v])<  \min \{ \ord (u), \ord (v)\} $ provided
$[u,v]\neq 0$;

(iv) $\ord (\s (u))=\ord (u)$ for all automorphisms $\s$ of the
Lie algebra $\ggu_n$ (Corollary \ref{d10Dec11}).


{\bf A classification of ideals of the Lie algebra $\ggu_n$}.
 The next lemma is the crucial  fact in the proof of Theorem
\ref{10Dec11}.

\begin{lemma}\label{a10Dec11}
Let $I$ be a nonzero ideal of the Lie algebra $\ggu_n$. Then
$I=\bigcup_{0\neq u\in I}I_{\ord (u)}$. In particular,
$(v)=I_{\ord (v)}$ for all nonzero elements  $v$ of $\ggu_n$ where
$(v)$ is the ideal of the Lie algebra $\ggu_n$ generated  by the
element $v$.
\end{lemma}

{\it Proof}. It suffices to show that the second statement holds,
that is $(v)=I_{\ord (v)}$,  since then
$$ I=\sum_{0\neq u\in I} (u) =\sum_{0\neq u\in I}I_{\ord (u)}=
\bigcup_{0\neq u\in I}I_{\ord (u)}.$$ To prove the equality
$(v)=I_{\ord (v)}$ we use the double induction: first on $n\geq 2$
and then, for a fixed $n$, the second induction on $\l = \ord
(v)$. Without loss of generality we may assume that the leading
coefficient of the element $v$ is $1$. Notice that $\ord (v)$ is
{\em not} a limit ordinal.

Let $n=2$. If $\l =1$ then $v=\der_2$, and so $(v) = K\der_2 =
I_1$, as required. Suppose that $\l >1$ and that the statement is
true for all non-limit ordinals $\l'$ such that $\l'<\l$.

Case 1: $\l = \ord (x_1^i\der_2)$ for some $i\geq 1$, i.e., $$v= x_1^i\der_2+\mu_{i-1}x_1^{i-1} \der_2+\cdots +\mu_0\der_2$$ for some
scalars $\mu_j\in K$. Let $\d = \ad (\der_1)$. Then $I_\l
\supseteq (v) \supseteq \sum_{j=0}^i K \d^j(v) = I_\l$, that is
$(v) = I_\l$.

Case 2: $\l = \ord (\der_1)=\o +1$, i.e.,  $v=\der_1+p\der_2$ for
some element $p\in P_1$. For all $i\geq 1$, $(v) \ni [ v,
x_1^i\der_2]=ix_1^{i-1}\der_2$. By Case 1, $\ggu_{2,1}\subseteq
(v)$, and so $\der_1\in (v)$. Therefore, $(v) = \ggu_2=I_{\o
+1}=I_{\ord (v)}$.

Suppose that $n>2$, and the statement holds for all $n'$ such that
$n'<n$. If $\l =1$ then $v\in K^*\der_n$, and so $(v) =K\der_n =
I_1$, as required. Suppose that $\l >1$ and that the statement is
true for all non-limit ordinals $\l'$ such that $\l '<\l$.

Case 1: $\l = \ord (x^\alpha \der_n)$ for some $0\neq \alpha \in
\N^{n-1}$, i.e.,  $v= x^\alpha \der_n+\cdots$ where the three dots
mean smaller terms. The next Claim follows from (\ref{Xai1})
and the definition of the well-ordering on  the canonical basis
$\CB_n$ of $\ggu_n$.

{\em Claim: For each  non-limit ordinal $\l'$ such that $\l'<\l$, there exist elements of $\CB_n$, say $a_1, \ldots , a_s$, of
type $X_{\beta , i}$ where $i\neq n$ such that }
$$\ord (\ad (a_1) \cdots \ad (a_s) (v))=\l'.$$

It follows from the Claim and the induction on $\l$  that $(v) =
I_{\ord (v)}$.

Case 2: $v=\der_i+\cdots $ for some $i$ such that $i<n$. For all
elements $x^\beta \der_{i+1}$ with $\beta \in \N^i$ and
$\beta_i\neq 0$,
$$(v)\ni [ v, x^\beta \der_{i+1}]=\beta_i x^{\beta
- e_i}\der_{i+1}+\cdots $$ and so $\ord ([ v, x^\beta
\der_{i+1}])=\ord ( x^{\beta -e_i} \der_{i+1})$. By induction,
$(v)\supseteq \cup_{\mu <\l } I_\mu=:J$. But $I_\l \supseteq
(v)=Kv+J= K\der_i +J=I_\l$, and so $(v) = I_\l$.

Case 3: $v=x^\alpha \der_i+\cdots$ for some $i$ such that $i<n$
and $\alpha \in \N^{i-1}\backslash \{ 0\}$. Notice that
$$ (v) \ni \prod_{j=1}^{i-1} \frac{\ad
(\der_j)^{\alpha_j}}{\alpha_j!}(v)= \der_i+\cdots . $$ By Case 2,
$I_{\ord (\der_i)} \subseteq (v)$. Since $\ggu_{n,i+1}\subseteq
I_{\ord (\der_i)}$, we have $\ggu_{n,i+1}\subseteq (v)$. By
(\ref{cuni1}), $\ggu_i\simeq \ggu_n/\ggu_{n,i+1}$. By the
assumption, $i<n$. By considering the element $v+\ggu_{n,i+1}\in
\ggu_i$, the  statement follows  by induction on $i$. $\Box $

Let $u,v\in \ggu_n\backslash \{ 0\}$. Then $(u)\subseteq (v)$ iff $\ord (u) \leq \ord (v)$.

\begin{theorem}\label{10Dec11}

\begin{enumerate}
\item The map (\ref{Jlid}) is a bijection. \item The Lie algebra
$\ggu_n$ is a uniserial, Artinian but not Noetherian Lie algebra
and its uniserial dimension is equal to $\udim (\ggu_n) = \ord
(\O_n) = \o^{n-1}+\o^{n-2}+\cdots +\o +1$.
\end{enumerate}
\end{theorem}

{\it Proof}. 1. Statement 1 follows at once from Lemma
\ref{a10Dec11}.

2. Statement 2 follows from statement 1. $\Box $


An ideal $\ga$ of a Lie algebra $\CG$ is called {\em proper}
(respectively, {\em co-finite}) if $\ga \neq 0, \CG$
(respectively, $\dim_K(\CG / \ga ) <\infty $).

\begin{corollary}\label{b10Dec11}

\begin{enumerate}
\item  The ideal $\ggu_{n,2}$ is the largest proper ideal of the
Lie algebra $\ggu_n$. \item The ideal $\ggu_{n,2}$ is the only
proper co-finite ideal of the Lie algebra $\ggu_n$, and
$\dim_K(\ggu_n/ \ggu_{n,2})=1$. \item The centre $Z(\ggu_n) =
K\der_n$ of the Lie algebra $\ggu_n$ is the least non-zero ideal
of the Lie algebra $\ggu_n$. \item The ideals $I_s:=
\sum_{i=0}^{s-1} Kx_1^i\der_n$ where $s=1,2,\ldots $ are the only
finite dimensional ideals of the Lie algebra $\ggu_n$, and
$\dim_K(I_s) = s$.
\end{enumerate}
\end{corollary}

{\it Proof}. All statements are easy corollaries of Theorem
\ref{10Dec11}.  $\Box $


{\bf The centralizers of the ideals of $\ggu_n$}.
In combination with Theorem \ref{10Dec11} the next proposition describes the centralizers of all the ideals of the Lie algebra $\ggu_n$. Notice that the centralizer of an ideal of Lie algebra is also an ideal.

\begin{proposition}\label{b24Jan12}
\begin{enumerate}
\item $$\Cen_{\ggu_n}(I_{\l , n})= \begin{cases}
K\der_n & \text{if }\l = \ord (\O_n) = \o^{n-1}+\cdots +\o +1,\\
I_{\o^i, n}=P_i\der_n& \text{if }\l\in (\o^{n-1}+\cdots + \o^{i+1},\o^{n-1}+\cdots + \o^i],  i=1, \ldots , n-2,\\
\oplus_{i=m+1}^n P_{i-1}\der_i & \text{if }\l \in (\o^{m-1}, \o^m], m=1, \ldots , n-1,\\
\ggu_n& \text{if }\l =1.\\
\end{cases}$$
\item The set $\CC (\ggu_n)$ of all the centralizers of the ideals of the Lie algebra $\ggu_n$ contains precisely $2n-1$ elements and the map $\CC (\ggu_n)\ra \CC (\ggu_n)$, $ C\mapsto \Cen_{\ggu_n}(C)$, is a bijection. Moreover, it is an inclusion reversing involution, i.e.,  $\Cen_{\ggu_n}(\Cen_{\ggu_n}(C))=C$ for all $C\in \CC (\ggu_n)$.
\item $\Cen_{\ggu_n}(\ggu_n) = K\der_n$, $\Cen_{\ggu_n}(K\der_n) = \ggu_n$, $\Cen_{\ggu_n}(P_{n-1}\der_n) = P_{n-1}\der_n$, $\Cen_{\ggu_n}(P_i\der_n) = \oplus_{j=i+1}^nP_{j-1}\der_j$ and $\Cen_{\ggu_n}( \oplus_{j=i+1}^nP_{j-1}\der_j)=P_i\der_n$ for $i=1, \ldots , n-2$.
\end{enumerate}
\end{proposition}

{\it Proof}. 1. Let $C_\l :=\Cen_{\ggu_n} (I_{\l , n})$. If $\l = \ord (\O_n )$ then $I_{\l , n} = \ggu_n$ and $C_\l = Z(\ggu_n) = K\der_n$, by Proposition \ref{a8Dec11}.(6). If $\l =1$ then $I_{1, n} = K\der_n = Z(\ggu_n)$, and so $C_1=\ggu_n$. We can assume that $\l \neq 1, \ord (\O_n)$. To prove the proposition we use induction on $n\geq 2$. For $n=2$, the only case to consider is when $\l \in (1, \o ]$. In this case,
$$I_{\l , 2}=\begin{cases}
\oplus_{i=0}^{\l -1} x_1^i \der_2& \text{if }1<\l <\o ,\\
P_1\der_2& \text{if }\l = \o , \\
\end{cases}
$$ and so $C_\l =P_1\der_2$.

Let $n>2$, and we assume that the proposition holds for all $n'<n$.

{\em Claim.} $C_{\o^m} = \oplus_{i=m+1}^nP_{i-1}\der_i$ for $m=1, \ldots , n-1$.

For $m=n-1$, $I_{\o^{n-1},n}=P_{n-1}\der_n$ and the claim is Lemma \ref{e10Dec11}.(2). Suppose that $m<n-1$. Then $I_{\o^m, n}=P_m\der_n\subset I_{\o^{n-1}, n}$, hence $C_{\o^m}\supseteq C_{\o^{n-1}}=P_{n-1}\der_n$. This means that the $\ggu_n$-module $I_{\o^m, n}$ is also a $\ggu_n/P_{n-1}\der_n$-module and a $\ggu_{n-1}$-module since $\ggu_{n-1}\simeq \ggu_n / P_{n-1}\der_n$. The map
$$I_{\o^m, n}=P_m\der_n\ra I_{\o^m, n-1}=P_m\der_{n-1}, \;\; p\der_n\mapsto p\der_{n-1},$$
 is an isomorphism of $\ggu_{n-1}$-modules (since $\der_n\in Z(\ggu_n)$ and $\der_{n-1}\in Z(\ggu_{n-1})$), and the claim follows by  induction on $n$.

Suppose that $1<\l \leq \o^{n-1}$, i.e.,  $\l \in (\o^{m-1}, \o^m]$, for some $m\in \{ 1, \ldots , n-1\}$. We have to show that $C_\l = C_{\o^m}$, see the claim.  By the claim, the inclusions $I_{\o^{m-1}, n}\subset I_{\l , n}\subseteq I_{\o^m , n}$ imply the inclusions  $C_{\o^{m-1}}\supseteq C_\l \supseteq C_{\o^m}$. Notice that $C_{\o^{m-1}}=C_{\o^m}\oplus P_{m-1}\der_m$, $x_m\der_n\in I_{\l , m}$ (since $\l \in (\o^{m-1}, \o^m]$), and, for all nonzero elements $p\in P_{m-1}$, $[p\der_m, x_m\der_n] = p\der_n\neq 0$. Hence, $C_\l = C_{\o^m}$, as required.

Suppose that $\o^{n-1}<\l <\ord (\O_n) = \o^{n-1}+\cdots +\o +1$, i.e.,  $\l \in (\o^{n-1}+\cdots + \o^{i+1},\o^{n-1}+\cdots + \o^i]$ for some $i\in \{ 1, \ldots , n-2\}$. Let $\l_i := \o^{n-1}+\cdots + \o^i$. Then $I_{\l , n}\subseteq I_{\l_i, n}=\oplus_{j=i+1}^nP_{j-1}\der_j=\Cen_{\ggu_n}(I_{\o^i, n})$, by the claim. Therefore, $C_\l \supseteq I_{\o^i , n}=P_i\der_n$. The inclusion $I_{\l , n}\supset P_{n-1}\der_n= I_{\o^{n-1}, n}$ implies the inclusion $C_\l \subseteq C_{\o^{n-1}}=P_{n-1}\der_n$, by the claim. Then the inclusion $\{ \der_{i+1}, \ldots , \der_n\} \subseteq I_{\l , n}$ implies that $C_\l \subseteq \Cen_{P_{n-1}}(\der_{i+1}, \ldots , \der_n)\der_n=P_i\der_n = I_{\o^i, n}$. This means that $C_\l = I_{\o^i, n}$.

2. Statement 2 follows from statement 3.

3. Statement 3 follows from statement 1. $\Box $


{\bf The central series and the central dimension}. For a Lie
algebra $\CG$ over a field $K$, let us define recursively its {\em
central series} $\{ Z^{(\l )} (\CG ) \}_{ \l \in \mW }$. Let    $Z^{(0
)} (\CG ):=Z(\CG )$.  If $\l $ is not a limit ordinal, that is $\l
= \mu +1$ for some $\mu \in \mW$,  then
$$Z^{(\l )} (\CG ):= \pi_\mu^{-1} (Z(\CG / Z^{(\mu )} (\CG )))\;\;{\rm
where }\;\; \pi_\mu : \CG \ra \CG / Z^{(\mu )} (\CG ), \;\;
a\mapsto a+Z^{(\mu )} (\CG ). $$ If $\l$ is a limit ordinal then
$Z^{(\l )} (\CG ):=\bigcup_{\mu <\l } Z^{(\mu )} (\CG )$. If $\l
\leq \mu$ then $Z^{(\l )} (\CG )\subseteq Z^{(\mu )} (\CG )$. So,
$\{ Z^{(\l )}(\CG )\}_{\l \in \mW}$ is an ascending class of
ideals of the Lie algebra $\CG$. Let $Z^{(\mW )} (\CG
):=\bigcup_{\l \in \mW} Z^{(\l )} (\CG )$.

$\noindent $

{\it Definition}. The minimal ordinal number $\l$ (if it exists)
such that $Z^{(\l )} (\CG )=Z^{(\mW )} (\CG )$ is called the {\em
central dimension} of the Lie algebra $\CG$ and is denoted by
$\cdim (\CG )$. In case, there is no such ordinal $\l $ we write
$\cdim (\CG ) =\mW$. The concept of the central dimension makes
sense for any class of algebras (associative, Jordan, etc.).

$\noindent $

The Lie algebra $\CG$ is {\em central}, that is $Z(\CG ) = 0$, iff
$\cdim (\CG ) =0$. So, the central dimension measures the deviation
from `being central.' The next theorem describes the central
series for the Lie algebra $\ggu_n$ and  gives $\cdim (\ggu_n)$.

\begin{theorem}\label{A10Dec11}
The central series $\{ Z^{(\l )} (\ggu_n) \}_{\l \in \mW}$
stabilizers precisely on $\ord (\O_n) = \o^{n-1}+\o^{n-2}+\cdots +
\o +1$ step, i.e.,  $\cdim (\ggu_n) = \ord (\O_n)$. Moreover, for
each $\l \in [1,\ord (\O_n)]$, $Z^{(\l )}(\ggu_n) = I_\l$. In particular, $Z^{(\cdim (\ggu_n))}(\ggu_n) = \ggu_n$.
\end{theorem}

{\it Proof}. It suffices to show that $Z^{(\l )}=I_\l$ for all $\l
\in [ 1,\ord (\O_n)]$ where $Z^{(\l )}:=Z^{(\l )}(\ggu_n)$. We use
an induction  on $\l$. The case $\l =1$ follows from Proposition
\ref{a8Dec11}.(6), $Z^{(1 )}=Z(\CG ) = K\der_n=I_1$. Suppose that
$\l >1$, and that the equality holds for all ordinals $\l'<\l$.

 If $\l$ is a limit ordinal then
 $$Z^{(\l )}=\bigcup_{\mu <\l } Z^{(\mu )}=\bigcup_{\mu <\l } I_\mu =I_\l.$$
 Suppose that $\l$ is not a limit ordinal, that is $\l = \mu +1$
 for some ordinal number $\mu$.

 Case 1: $\l \leq \o^{n-1}$. Then necessarily $\l <\o^{n-1}$ since
 $\l $ is not a limit ordinal but $\o^{n-1}$ is. Clearly, $I_\l
 \subseteq Z^{(\l )}$, by (\ref{Xai1}) (since $Z^{(\mu )}=I_\mu$, by induction).  Suppose that $I_\l \neq Z^{(\l
 )}$ then necessarily $I_{\l +1} \subseteq Z^{(\l )}$ (by Theorem
 \ref{10Dec11}.(1)). If $\l = \Xai$ for some $\alpha \in \N^{n-1}$
 (recall that we identified $(\O_n , \leq )$ and $(\CB_n , \leq )$)
 then $\l +1= X_{\alpha +e_1, n}$ and $ X_{\alpha +e_1, n}\in
 I_{\l +1}$, but $[\der_1,  X_{\alpha +e_1,n}]=(\alpha_1+1)  X_{\alpha ,
 n}\not\in Z^{(\mu )}$, a contradiction. Therefore, $Z^{(\l
 )}=I_\l$.

Case 2: $\l > \o^{n-1}$. Notice that  $Z^{(\l )}\supseteq I_\l
\supsetneqq I_{\o^{n-1}}= \ggu_{n,n}$ and $\ggu_n /
\ggu_{n,n}\simeq \ggu_{n-1}$, by (\ref{cuni1}). Now, we  complete
the argument by induction on $n$. The base of the induction $n=2$
is covered by the Cases 1 and 2 for $n=2$ above.
 $\Box $


An ideal $I$  of a Lie algebra $\CG$ is called a {\em
characteristic ideal} if it is invariant under all the
automorphisms of the Lie algebra $\CG$, that is $\s (I)=I$ for all
$\s \in \Aut_K(\CG )$. It is obvious that an ideal $I$ is a
characteristic ideal iff $\s (I) \subseteq I$ for all $\s \in
\Aut_K(\CG )$.

\begin{corollary}\label{c10Dec11}
All the ideals of the Lie algebra $\ggu_n$ are characteristic
ideals.
\end{corollary}

{\it Proof}. By the very definition, all the ideals in the central
series $\{ Z^{(\l )}(\ggu_n)\}$ are characteristic ideals of
$\ggu_n$ but these are all the ideals of the Lie algebra $\ggu_n$,
by Theorem \ref{A10Dec11}. We can also deduce the corollary from Theorem \ref{10Dec11}.(2).  $\Box $


\begin{corollary}\label{d10Dec11}
For all nonzero elements $u\in \ggu_n$ and all automorphisms $\s $
of the Lie algebra $\ggu_n$, $\ord (\s (u)) = \ord (u)$.
\end{corollary}

{\it Proof}. Notice that 
\begin{equation}\label{ordul}
\ord (u) = \min \{ \l \in [ 1, \ord (\O_n) ] \, | \, u\in I_\l \}.
\end{equation}
Now, the statement follows from the fact that all the ideals
$I_\l$ are characteristic (Corollary \ref{c10Dec11}). $\Box $


{\bf The subalgebra $U_n'$ of the Weyl algebra $A_n$}.
Let be $\ga $ an ideal of a Lie algebra $\CG$. Then the ideal $U(\CG ) \ga U(\CG )$ of the universal enveloping algebra $U(\CG )$ generated by $\ga$ is equal to $U(\CG ) \ga = \ga U(\CG )$. The chain of ideals of the Lie algebra $\ggu_n$
$$ I_1\subset \cdots \subset I_\l \subset\cdots \subset \ggu_n = I_{\ord (\O_n )}, \;\; \l \in [1, \ord (\O_n)], $$ yields
the chain of subalgebras and the chain of ideals of the algebra $U_n$,  respectively,
$$ U( I_1)\subset \cdots \subset U( I_\l )\subset \cdots \subset U_n \;\;\;\; {\rm and}\;\;\;\;
 U_nI_1\subset \cdots \subset U_nI_\l\subset \cdots \subset U_n.$$
Let $U_{n,\l }':=\chi_n(U(I_\l ))$ and $I_\l':=\chi_n(U_nI_\l )$ where $\chi_n$ is the algebra homomorphism (\ref{UnAn}). Then
\begin{equation}\label{UpnLp}
U_{n,1 }'\subseteq \cdots \subseteq U_{n,\l }'\subseteq \cdots \subseteq U_n'=\chi_n (U_n)
\end{equation}
is a chain of subalgebras of the algebra $U_n'$, and
\begin{equation}\label{UpnLp1}
I_{n,1 }'\subseteq \cdots \subseteq I_{n,\l }'\subseteq \cdots \subseteq U_n'
\end{equation}
is a chain of ideals of the algebra $U_n'$.

\begin{proposition}\label{a14Jan12}
\begin{enumerate}
\item All the inclusions in (\ref{UpnLp}) are strict inclusions. The algebra $U_n'$ is not a finitely generated algebra.
\item All the inclusions in (\ref{UpnLp1}) are strict inclusions. In particular, the algebra $U_n'$ is neither left nor right  Noetherian  algebra and  does not satisfy the ascending chain condition on ideals.
\end{enumerate}
\end{proposition}

{\it Proof}. 1. We use induction on $n$. Let $n=2$. For each $i\in [1, \o )$, the algebra $U_{2,i}'$ is generated by the commuting elements $x_1^j\der_2$, where $0\leq j\leq i-1$ since $$ x_1^j\der_2\cdot x_1^k\der_2=x_1^{j+k}\der_2^2\;\;\; {\rm for \; all}\;\; j,k.$$ These equalities mean that the algebra $U_{2,i}'$ is isomorphic to the monoid algebra $K\CM_i$ (via $x_1^j\der_2\mapsto  (j,1)$) where $\CM_i$ is the submonoid of $(\N^2, +)$ generated by the elements $(j,1)$ where $0\leq j\leq i-1$. It follows that
$$ x_1^{i-1}\der_2\in U_{2,i}'\backslash U_{2,i-1}'\;\; {\rm for \; all}\;\; i=2,3,\ldots .$$
This means that the inclusions $U_{2,1}'\subset U_{2,2}' \subset\cdots \subset U_{2,i}'\subset \cdots \subset U_{2, \o}'=\cup_{i\geq 1}U_{2,i}'$ are strict inclusions and so the algebra $U_{2, \o}'$ is not a finitely generated algebra. Since $\der_1\in U_2'\backslash U_{2, \o }'$, $U_2'=U_{2, \o}'[\der_1; \ad (\der_1)]$ is the skew polynomial algebra and $[\der_1, U_{2,i}']\subseteq U_{2,i}'$ for all $i\geq 1$,  statement 1 is true for $n=2$.

Suppose that $n>2$  and we assume that statement 1 is true for all $n'<n$. Let $ \l \in [1, \o^{n-1})$. Then $\l +1$ is not  a limit ordinal. Notice that $[1, \o^{n-1})\subseteq [ 1, \ord (\O_n))$. Hence, $\l+1=(\alpha , n)$. The elements $\{ x^{\beta}\der_n\, | \, \beta \in \N^{n-1}\}$ commute. Moreover,
$$ x^{\beta} \der_n \cdot x^\g \der_n = x^{\beta +\g } \der_n^2\;\;\; {\rm for \; }\;\; \beta , \g \in \N^{n-1}.$$
Therefore, the algebra $U_{n,\l}'$ is a commutative algebra which is isomorphic to the monoid algebra $K\CM_{n, \l}$ (via $x^{\beta} \der_n\mapsto (\beta , 1)$) where $\CM_{n, \l}$ is the submonoid of $(\N^n, +)$ generated by the elements $\{ (\beta, 1)\, | \, (\beta , n)\leq \l \}$. It follows that
$$x^\alpha\der_n \in U_{n,\l +1}'\backslash U_{n, \l }'\;\;\; {\rm for \; all}\;\; \l \in [1, \o^{n-1}).$$ This means that (the inclusions are strict) $$U_{n,1}'\subset U_{n,2}'\subset \cdots \subset U_{n,\l }'\subset \cdots \subset U_{n,\o^{n-1} }'=\bigcup_{\l \in [1,\o^{n-1})}U_{n, \l}'.$$
Let $\l \geq \o^{n-1}$. By Theorem \ref{12Jan12}.(1) and (\ref{Jlid}),
\begin{equation}\label{Upnlp2}
U_n'/I_{n,\o^{n-1}}'\simeq U_{n-1}', \;\;
(U_{n,\o^{n-1}+\mu}'+I_{n,\o^{n-1}}')/I_{n,\o^{n-1}}'\simeq U_{n-1, \mu}'\;\;\; {\rm for \; all}\;\; \mu\in[1,\ord (\O_{n-1})].
\end{equation}
By induction on $n$, statement 1 holds.

2. We use freely the facts proved in the proof of statement 1.  We use induction on $n$. Let $n=2$. For each $i\in [1,\o )$, the ideal $I_{2,i}'$ of the algebra $U_2'$ is generated by the commuting elements $x_1^j\der_2$ where $0\leq j\leq i-1$. By Corollary \ref{x12Jan12}, the set $W_2'=\{ \der_1^i\der_2^j, \der_1^ix_1^k\der_2\cdot \der_2^j\, | \, i,j\in \N, k\geq 1\}$ is a $K$-basis for the algebra $U_2'$. For each $i\in [1,\o )$,
$$ I_{2,i}'=\sum_{j=0}^{i-1}U_2'x_1^j\der_2=\bigoplus_{i\geq 0}\bigoplus_{k\geq 0, j\geq 2}K\der_1^ix_1^k\der_2^j\oplus \bigoplus_{i\geq 0}\bigoplus_{k=0}^{i-1}K\der_1^ix_1^k\der_2.$$ It follows that
$$x_1^{i-1}\der_2\in I_{2,i}'\backslash I_{2,i-1}', \;\; i=2,3,\ldots .$$
This means that
$$I_{2,1}'\subset I_{2,2}'\subset \cdots \subset I_{2,i}'\subset\cdots \subset I_{2,\o}'=\bigcup_{i\geq 1}I_{2,i}'.$$
 Since $1\in U_2'\backslash I_{2,\o +1}'$ and $\der_1\in I_{2, \o+1}'\backslash I_{2, \o}'$, statement 2 is true for $n=2$.

Suppose that $n>2$  and we assume that statement 2 is true for all $n'<n$. Let $ \l \in [1, \o^{n-1})$. Then $\l +1$ is not  a limit ordinal. Notice that $[1, \o^{n-1})\subseteq [ 1, \ord (\O_n))$. Hence, $\l+1=(\alpha , n)$.  It follows from the equality $I_{n,\l +1}'=U_n'\chi_n(I_{\l +1})=\sum_{(\beta , n)\leq \l +1}U_n'x^\beta \der_n$ and (\ref{Xai1}) that
\begin{equation}\label{Wpn1}
I_{n,\l +1}'\cap(\bigoplus_{\beta \in \N^{n-1}}K x^{\beta}\der_n)=\bigoplus_{(\beta , n)\leq \l +1}Kx^{\beta}\der_n.
\end{equation}
It follows that
$$x^\alpha\der_n \in I_{n,\l +1}'\backslash I_{n, \l }'\;\;\; {\rm for \; all}\;\; \l \in [1, \o^{n-1}).$$ This means that $$I_{n,1}'\subset I_{n,2}'\subset \cdots \subset I_{n,\l }'\subset \cdots \subset I_{n,\o^{n-1} }'=\bigcup_{\l \in [1,\o^{n-1})}I_{n, \l}'.$$
Let $\l >\o^{n-1}$. In view of (\ref{Upnlp2}), we also have
\begin{equation}\label{Upnlp3}
I_{n,\o^{n-1}+\mu}'/I_{n,\o^{n-1}}'\simeq I_{n-1, \mu}'\;\;\; {\rm for \; all}\;\; \mu\in[1,\ord (\O_{n-1})].
\end{equation}
By induction on $n$, statement 2 holds. $\Box $


{\bf The Heisenberg Lie subalgebras of $\ggu_n$}.
 Let $\gl_n(K)=\bigoplus_{i,j=1}^nKE_{ij}$ be the $n\times n$
matrix Lie algebra where $E_{ij}$ are the matrix units. For each
$n\geq 2$, let $\UT_n(K)=\bigoplus_{i\leq j}KE_{ij}$ be  the {\em
upper triangular} $n\times n$ matrix Lie algebra,
$\UT_n(K)\subseteq \gl_n(K)$. The $K$-linear map
$$\UT_n(K)\ra
\ggu_n, \;\; E_{ij}\mapsto x_i\der_j,$$ is a Lie algebra monomorphism.
We identify the Lie algebra $\UT_n(K)$ with its image in $\ggu_n$.

The {\em Heisenberg Lie algebra} $\CH_n$ is a $2n+1$ dimensional
Lie algebra with a $K$-basis $X_1, \ldots , X_n$, $ Y_1, \ldots ,
Y_n, Z$ where $Z$ is a central element of $\CH_n$,
$$ [Y_i, X_j]=\d_{ij}Z, \;\; [X_i,X_j]=0, \;\; [Y_i,Y_j]=0, \;\;
i,j=1, \ldots , n$$ where $\d_{ij}$ is the Kronecker
delta-function.  The $K$-linear map
$$\CH_{n-1}\ra \ggu_n , \;\; X_i\mapsto x_i\der_n, \;\; Y_i\mapsto
\der_i, \;\; Z\mapsto \der_n,\;\; i=1, \ldots , n-1$$ is a Lie
algebra monomorphism. We identify the Heisenberg Lie algebra $\CH_{n-1}$ with its image in
$\ggu_n$. Since $\CH_{i-1}\subseteq \ggu_i\subseteq \ggu_n$ for $i=2, \ldots , n-1$, the Lie algebra $\ggu_n$ contains the Heisenberg Lie algebras $\CH_i$, $i=1, \ldots , n-1$. For all natural integers $i$ and $j$ with $i\neq j$, $\CH_i\cap \CH_j = \sum_{k=1}^{\min \{ i,j\}}K\der_k$.

Let $A$ be an algebra or a Lie algebra and $M$ be its module. Then $\ann_A(M):=\{ a\in A \, | \, aM=0\}$ is called the {\em annihilator} of the $A$-module $M$. It is an ideal of $A$. A module is called {\em faithful} if its annihilator is 0.

{\bf The $\ggu_n$-module $P_n$}. To say that $P_n$ is a $\ggu_n$-module is the same as to say that $P_n$ is a $U_n'$-module. The importance of this obvious observation is that we can use the relations of the Weyl algebra $A_n$ in various computations with $U_n'$  as $U_n'\subseteq A_n$.
For each $n\geq 2$, $\ggu_n$ is a Lie subalgebra of the Lie algebra $\ggu_{n+1} = \ggu_n\oplus P_n\der_{n+1}$, $P_n\der_{n+1}$ is an ideal of the Lie algebra $\ggu_{n+1}$ and $[P_n\der_{n+1},P_n\der_{n+1}]=0$. In particular, $P_n\der_{n+1}$ is a left $\ggu_n$-module where the action of the Lie algebra $\ggu_n$ on $P_n\der_{n+1}$ is given by the rule: $uv:= [u,v]$ for all $u\in \ggu_n$ and $v\in P_n\der_{n+1}$. Recall that the polynomial algebra $P_n$ is a left $\ggu_n$-module.

\begin{lemma}\label{aa21Jan12}
\begin{enumerate}
\item The $K$-linear map $P_n\ra P_n\der_{n+1}$, $p\mapsto p \der_{n+1}$, is a $\ggu_n$-module isomorphism.
\item The $\ggu_n$-module $P_n$ is an indecomposable, uniserial $\ggu_n$-module, $\udim (P_n)=\o^n$ and $\ann_{\ggu_n}(P_n)=0$.
\item The set $\{ P_{\l , n}:=\oplus_{\alpha \in \N^n, (\alpha , n+1)\leq \l} Kx^\alpha \, | \, \l \in [ 1, \o^n]\}$ is the set of all the nonzero  $\ggu_n$-submodules of $P_n$; $P_{\l , n}\subset P_{\mu , n}$ iff $\l < \mu$; $\udim (P_{\l , n})=\l$ for all $\l \in [1, \o^n]$.
    \item All the $\ggu_n$-submodules of $P_n$ are pairwise non-isomorphic, indecomposable, uniserial $\ggu_n$-modules.
\end{enumerate}
\end{lemma}

{\it Proof}. 1. It is obvious that the map is a bijection and a $\ggu_n$-homomorphism.

3.  By statement 1, the $\ggu_n$-module $P_n$ can be identified with the ideal $P_n\der_{n+1}$ of the Lie algebra $\ggu_{n+1}$. Under this identification every $\ggu_n$-submodule of $P_n$ becomes an ideal of the Lie algebra $\ggu_{n+1}$ in $P_n\der_{n+1}$, and vice versa.  Now, statement 3 follows from statement 1 and the classification of ideals of the Lie algebra $\ggu_{n+1}$ (Theorem \ref{10Dec11}.(1)).

2. By statement 3, the $\ggu_n$-module $P_n$ is uniserial, hence indecomposable, and  $\udim (P_n)=\o^n$. The Weyl algebra $A_n$ is a simple algebra, hence $\ann_{A_n}(P_n)=0$. Then $\ann_{\ggu_n}(P_n)=\ggu_n\cap \ann_{A_n}(P_n)=0$ since $\ggu_n\subseteq A_n$.

4. By statement 3, $\udim (P_{\l , n})=\l$. Hence $P_{\l , n}\not\simeq P_{\mu , n}$ for all $\l \neq \mu$. The rest is obvious (see statement 2). $\Box $


The next corollary describes the annihilators of all the $\ggu_n$-submodules of $P_n$. In particular, it classifies the faithful $\ggu_n$-submodules of the $\ggu_n$-module $P_n$.

\begin{corollary}\label{a24Jan12}
\begin{enumerate}
\item The $\ggu_n$-submodule $P_{\l ,  n}$ of $P_n$ is a faithful submodule iff $\l \in [\o^{n-1}+1, \o^n]$.
\item $\ann_{\ggu_n}(P_{\o^{n-1}, n}) = P_{n-1}\der_n$.
\item $\ann_{\ggu_n}(P_{\l , n}) =\begin{cases}
0& \text{if }\l \in (\o^{n-1}, \o^n],\\
\oplus_{i=m+1}^nP_{i-1}\der_i& \text{if }\l \in (\o^{m-1}, \o^m], m=1, \ldots , n-1, \\
\ggu_n & \text{if }\l =1.\\
\end{cases}$
\end{enumerate}
\end{corollary}

{\it Proof}. 2. The inclusion ``$\supseteq $'' is obvious. The reverse inclusion  ``$\subseteq $'' follows from the facts that $P_{n-1}$ is a faithful $\ggu_{n-1}$-module (Lemma \ref{aa21Jan12}.(2)), $\ggu_{n-1}\simeq \ggu_n / P_{n-1}\der_n$ and $P_{n-1}=P_{\o^{n-1}, n}$ is an $\ggu_n/P_{n-1}\der_n$-module.

1. Notice that if $N$ is  a submodule of $M$ then $\ann (N)\supseteq \ann (M)$. In view of this fact and statement 2, to finish the proof of statement 1 it suffices to show that $\ggu_n$-module $P_{\o^{n-1}+1, n}$ is faithful. Since $P_{\o^{n-1}, n}\subseteq P_{\o^{n-1}+1, n}$, we have $\ga :  = \ann_{\ggu_n} (P_{\o^{n-1}+1, n})\subseteq  \ann_{\ggu_n} (P_{\o^{n-1}, n})=P_{n-1}\der_n$ Since $x_n\in P_{\o^{n-1}+1, n}$ and, for all nonzero  elements $p\in P_{n-1}$, $p\der_n *(x_n) = p\neq 0$, we have $\ga =0$.

3. We use induction on $n\geq 2$. The initial case is $n=2$, and  there are three cases to consider: $\l \in (\o ,  \o^2]$, $\l = (1, \o ]$ and $\l =1$. The first case is statements 1. The last case is obvious.  In the second case,
$$\ann_{\ggu_2} (P_{m,2}) = \ann_{\ggu_2} (\oplus_{i=0}^{m-1}Kx_1^i) = K[x_1]\der_2.$$ Let $n>2$, and  we assume that statement 3 holds for all $n'<n$. If $\l \in (\o^{n-1}, \o^n]$ then, by statement 1, the $\ggu_n$-module $P_{\l , n}$ is faithful. If $\l =1$ then $\ann_{\ggu_n} (P_{1, n})=\ann_{\ggu_n} (K)=\ggu_n$. If $1<\l \leq \o^{n-1}$, i.e.,  $\l \in (\o^{m-1}, \o^m]$ for some $m\in \{ 1,\ldots , n-1\}$, then $P_{\l , n} \subseteq P_{n-1}$. By statement 2, $\ga := \ann_{\ggu_n}(P_{\l , n})\supseteq \ann_{\ggu_n}(P_{\o^{n-1} , n})=P_{n-1}\der_n$. Since $\ggu_{n-1}\simeq \ggu_n / P_{n-1}\der_n$, the inclusion $P_{\l , n} \subseteq P_{n-1}$ is an inclusion of $\ggu_{n-1}$-modules. Moreover, the $\ggu_{n-1}$-submodule $P_{\l , n}$ of $P_{n-1}$ can be identified with the $\ggu_{n-1}$-submodule $P_{\l , n-1}$ of $P_{n-1}$. Now, the result follows by induction on $n$.   $\Box $


\begin{corollary}\label{xa24Jan12}
Let $n\geq 3$. Then the ideal $I_{\o^{n-2}+1}=Kx_{n-1}\der_n+\sum_{\alpha \in \N^{n-2}}Kx^\alpha \der_n$ is the least ideal $I$ of the Lie algebra $\ggu_n$ which is a faithful $\ggu_{n-1}$-module (i.e.,   $\ann_{\ggu_{n-1}}(I)=0$).
\end{corollary}

{\it Proof}. By Lemma \ref{aa21Jan12}.(1) and Corollary \ref{a24Jan12}.(1), the $\ggu_{n-1}$-module $I_{\o^{n-2}+1}$ is faithful but its predecessor $I_{\o^{n-2}}$ is not as  $\ann_{\ggu_{n-1}}(I_{\o^{n-2}})=\ann_{\ggu_{n-1}}(\sum_{\alpha \in \N^{n-2}}Kx^\alpha \der_n) \ni \der_{n-1}$. $\Box $


The inclusion of Lie algebras $\ggu_n\subseteq \ggu_{n+1}=\ggu_{n}\oplus P_{n-1}\der_n$ respects the total orderings on the bases $\CB_{n}$ and $\CB_{n+1}$. The $\ggu_n$-module isomorphism $P_n\ra P_n\der_{n+1}$, $p\ra p\der_{n+1}$, (Lemma \ref{aa21Jan12}.(1)), induces the total ordering on the monomials $\{ x^\alpha \}_{\alpha \in \N^n }$ of the polynomial algebra $P_n$ by the rule $x^\alpha >x^\beta $ iff $X_{\alpha , n+1}>X_{\beta , n+1}$ iff $\alpha_n = \beta_n$, $\alpha_{n-1} = \beta_{n-1}, \ldots , \alpha_{m+1} = \beta_{m+1}$ and $\alpha_m > \beta_m$ for some $m\in \{ 1, \ldots , n\}$.  This is the, so-called, {\em reverse lexicographic ordering} on $\{ x^\alpha \}_{\alpha \in \N^n}$ or $\N^n$ ($\alpha >\beta $ iff $x^\alpha > x^\beta $). By Lemma \ref{a18Dec11}.(3,4) and Lemma \ref{aa21Jan12}.(1), if $x^\alpha > x^\beta$ (where $\alpha , \beta \in \N^n$) then

(i) $\der_k*x^\alpha > \der_k*x^\beta$ for all $k=1, \ldots , n$ such that $\alpha_k\neq 0$; and

(ii) $X_{\g , k} *x^\alpha >X_{\g , k} *x^\beta$  for all  $k=1, \ldots , n-1$ and $\g \in \N^{k-1}$  such that $X_{\g , k} *x^\alpha \neq 0$, i.e.,  $\alpha_k\neq 0$.

Let $\CM (P_n)$ be the set of all nonzero submodules of the $\ggu_n$-module $P_n$, it is a well-ordered set with respect to $\subseteq$. By Lemma \ref{aa21Jan12}.(3), the map
\begin{equation}\label{xJlid}
\kappa_n :[1, \o^n]\ra \CM (P_n), \;\; \l \mapsto P_{\l , n},
\end{equation}
is the isomorphism of well-ordered sets. Each nonzero polynomial $p\in P_n$ is the unique sum $$\l_\alpha x^\alpha +\sum_{\alpha > \beta} \l_\beta x^\beta\;\; {\rm  where}\;\; \l_\alpha \in K^*\;\; {\rm  and}\;\; \l_\beta \in K.$$ The elements $\l_\alpha x^\alpha$ and $\l_\alpha$ are called the {\em leading term} and the {\em leading coefficient}  of the polynomial $p$ respectively (with respect to the well-ordering $>$). The ordinal number $$\ord ( p) := \ord ([1, \alpha ])$$ (where $[1, \alpha ] \subseteq [ 1, \o^n]$) is called the {\em ordinal degree} of $p$. For all nonzero polynomials $p,q\in P_n$ and $\l \in K^*$,

(i) $\ord (p+q) \leq \max \{ \ord (p) , \ord (q)\}$ provided $p+q\neq 0$;

(ii) $\ord ( \l p) = \ord (p)$; and

(iii) $\ord (u*p)<\ord (p)$ for all $u\in \ggu_n$ such that $u*p\neq 0$.

For each ordinal $\l \in [ 1,\o^n)$, there is the unique presentation $$\l = \alpha_n\o^{n-1}+\alpha_{n-1}\o^{n-2}+\cdots +\alpha_2\o +\alpha_1$$ where $\alpha_i\in \N$ and not all $\alpha_j$ are equal to zero (notice the shift by 1 of the indices in the coefficients $\alpha_i$). Then 
\begin{equation}\label{Pln1}
P_{\l , n}=\begin{cases}
\oplus \{ Kx^\beta \, | \, \beta \in \N^n, x^\beta \leq x_1^{-1}x^\alpha \} & \text{if }\alpha_1\neq 0,\\
\oplus \{ Kx^\beta \, | \, \beta \in \N^n, x^\beta < x^\alpha \}& \text{if }\alpha_1= 0,\\
\end{cases}
\end{equation}
where $x^\alpha =\prod_{i=1}^n x_i^{\alpha_i}$. Notice that $\alpha_1\neq 0$ iff $\l $ is not a limit ordinal. The vector space $P_{\l , n}$ has the largest monomial  iff $\l$ is not a limit ordinal, and in this case $x_1^{-1}x^\alpha$ is the largest monomial of $P_{\l , n}$.
  The ordinal number $\l \in [ 1,\o^n)$  is the unique sum  $$\l = \alpha_m\o^{m-1}+\alpha_{m-1}\o^{m-2}+\cdots +\alpha_j\o^{j-1}, \;\; \alpha_m\neq 0 , \;\; \alpha_j\neq 0,  $$ where $\alpha_i\in \N$, $1\leq m \leq n$ and $1\leq j \leq m$.  The positive integers $\alpha_m$ and $\alpha_j$ are called the {\em multiplicity} and the {\em co-multiplicity} of the ordinal number $\l$. The natural numbers  $m-1$ and $j-1$ are called the {\em degree} and the {\em co-degree} of the ordinal number $\l$.

\begin{lemma}\label{a26Feb12}
Let $\l \in [1,\o^n)$, i.e.,  $\l = \alpha_m\o^{m-1}+\alpha_{m-1}\o^{m-2}+\cdots + \alpha_j\o^{j-1}$ with $\alpha_i\in \N$, $\alpha_m\neq 0$, $\alpha_j\neq 0$, $m\leq n$ and $j\geq 1$. Then
\begin{eqnarray*}
 P_{\l , n}&=&\sum_{i=0}^{\alpha_m-1}x_m^iP_{m-1}+x_m^{\alpha_m}\sum_{0\leq i \leq \alpha_{m-1}-1}x_{m-1}^iP_{m-2}+\cdots \\
&+&x_m^{\alpha_m} x_{m-1}^{\alpha_{m-1}} \cdots x_{k+1}^{\alpha_{k+1}}
\sum_{0\leq i \leq \alpha_k-1}x_k^iP_{k-1}+\cdots +
x_m^{\alpha_m} x_{m-1}^{\alpha_{m-1}} \cdots x_{j+1}^{\alpha_{j+1}}
\sum_{0\leq i \leq \alpha_j-1}x_j^iP_{j-1}.
 \end{eqnarray*}
 In particular, $P_{\alpha_m\o^{m-1}, n}=\sum_{i=0}^{\alpha_m-1}x_m^iP_{m-1}$.
\end{lemma}

{\em Remark}. If $\alpha_k=0$ then the corresponding summand is absent.

{\it Proof}. Trivial.  $\Box $


Let $V$ be a vector space over the field $K$, a linear map $\v :
V\ra V$ is called a {\em Fredholm map/operator} if it has finite
dimensional kernel and cokernel, and then $$\ind (\v ) := \dim_K(
\ker (\v )) - \dim_K(\coker (\v ))$$ is called the {\em index} of
the map $\v$. Let $\CF (V)$ be the set of all Fredholm linear maps
in
$V$. In fact,  it  is a monoid since 
\begin{equation}\label{indpps}
\ind (\v \psi ) = \ind (\v ) +\ind (\psi )\;\; {\rm for \; all}
\;\; \v , \psi \in \CF (V).
\end{equation}
Let $V$ be a vector space with a countable $K$-basis $\{ e_i\}_{i\in \N}$ and $\der$ be a $K$-linear map on $V$ given by the rule $\der e_i=e_{i-1}$ for all $i\in \N$ where $e_{-1}:=0$.  For example, $V=K[x]$, $e_i:=\frac{x^i}{i!}$  and $\der =\frac{d}{dx}$. The subalgebra of $\End_K(V)$ generated by the map $\der $ is a polynomial algebra $K[\der ]$. Since the map $\der$ is a locally nilpotent map the algebra $\End_K(V)$ contains the algebra  $K[[\der]]=\{ \sum_{i=0}^\infty \l_i \der^i\, | \, \l_i\in K\}$ of formal power series in $\der$. The set $K[[\der]]^*=\{ \sum_{i=0}^\infty \l_i \der^i\in K[[\der ]]\, | \, \l_0\in K^*\}$ is the group of units of the algebra $K[[\der ]]$. The vector space $V$ is a $K[\der ]$-module.

\begin{lemma}\label{a26Jan12}

\begin{enumerate}
\item $\End_{K[\der ]}(V)=\{ \v \in \End_K(V), \v (e_i) = \sum_{j=0}^i\l_je_{i-j}=(\sum_{j\geq 0}\l_j\der^j)(e_i), i\in \N \, | \, \l_i\in K, i\in \N\}=K[[\der ]]$  and
    $\Aut_{K[\der ]}(V)=\{ \v \in \End_K(V), \v (e_i) = \sum_{j=0}^i\l_je_{i-j}, i\in \N \, | \, \l_0\neq 0,  \l_i\in K, i\in \N\}=K[[\der ]]^*$.
\item In particular, all nonzero elements of $\End_{K[\der ]}(V)$ are surjective Fredholm maps.
\item Let $a=\sum_{i\geq d} \l_i \der^i \in K[[\der ]]$ and $\l_d\neq 0$. Then $d=\ind (a) = \dim_K(\ker (a))$.
\end{enumerate}
\end{lemma}

{\it Proof}. Trivial.  $\Box $


Let $V$ be a vector space with a $K$-basis $\{ e_i\}_{i\in \N}$ and $x$ be a $K$-linear map on $V$ defined as follows $xe_i:=e_{i+1}$  for all $i\in \N$. For example, $V=K[x]$, $e_i:= x^i$  and $x : K[x]\ra K[x]$, $p\mapsto xp$. The polynomial algebra $K[x]$ is a subalgebra of the algebra $\End_K(V)$.

\begin{lemma}\label{a27Jan12}

\begin{enumerate}
\item $\End_{K[x ]}(V)=\{ p:V\ra V, a\mapsto pa\, | \, p\in K[x]\} \simeq K[x]$ (via $p\mapsto p$). In particular, all nonzero elements of $\End_{K[x]}(V)$ are injective maps.
\item  $\Aut_{K[x]}(V) =K^*$.
\item  For all $p\in \End_{K[x]}(V) = K[x]$, $\deg (p) = -\ind (p) = \dim_K(\coker (p))$.
\end{enumerate}
\end{lemma}

{\it Proof}. Trivial.  $\Box $


The next proposition describes the algebra of all the $\ggu_n$-homomorphisms (and its group of units) of the $\ggu_n$-module $P_n$. Clearly, $K[x_n]\subseteq P_n$ is the inclusion of $K[\der_n]$-modules and $K[[\frac{d}{dx_n}]]= \End_{K[\der_n]}(K[x_n])$ (Lemma \ref{a26Jan12}.(1)). The $K$-derivation $\frac{d}{dx_n}$ of the polynomial algebra $K[x_n]$ is also denoted by $\der_n$.

\begin{proposition}\label{A26Jan12}

\begin{enumerate}
\item The map $\End_{\ggu_n}(P_n) \ra \End_{K[\der_n]}(K[x_n]) = K[[\frac{d}{d x_n}]]$, $\v \mapsto \v |_{K[x_n]}$, is a $K$-algebra isomorphism with the inverse map $\v'\mapsto \v$ where $\v (x^\beta x_n^i):= X_{\beta , n} \v'(\frac{x_n^{i+1}}{i+1})$ for all $\beta \in \N^{n-1}$ and $ i\in \N$.
\item The map $\Aut_{\ggu_n} (P_n) \ra \Aut_{K[\der_n]}(K[x_n])= K[[\frac{d}{dx_n}]]^*$, $\v \mapsto \v |_{K[x_n]}$, is a group  isomorphism with the inverse map as in statement 1.
\item Every nonzero map $\v \in \End_{\ggu_n} (P_n)$ is a surjective map with kernel $\ker (\v ) =\oplus_{i=0}^{d-1}P_{n-1}x_n^i$ where $d=\ind (\v |_{K[x_n]}) = \dim_K(\ker (\v |_{K[x_n]}))$.
    \item For all natural numbers $d\geq 1$, $(\frac{\der}{\der x_n})^d\in \End_{\ggu_n}(P_n)$ and $ \ker_{P_n} (\frac{\der}{\der x_n})^d=\oplus_{i=0}^{d-1} P_{n-1}x_n^i$. In particular, $(\frac{\der}{\der x_n})^d: P_n/\oplus_{i=0}^{d-1} P_{n-1}x_n^i\ra P_n$ is a $\ggu_n$-module isomorphism.
\end{enumerate}
\end{proposition}

{\it Proof}. 1. By Lemma \ref{a26Jan12}, $\End_{K[\der_n]}(K[x_n])=K[[\frac{d}{dx_n}]]$. Since $K[x_n]=\cap_{i=1}^{n-1}\ker_{P_n}(\der_i)$ and the maps $\der_1, \ldots , \der_n$ commute, the restriction map $\v \mapsto \v |_{K[x_n]}$ is a well defined  $K$-algebra homomorphism. For all elements $\beta \in \N^{n-1}$ and $i\in \N$,
$$\v (x^\beta x_n^i) = \v ( X_{\beta , n}* \frac{x_n^{i+1}}{i+1}) = X_{\beta , n} \v (\frac{x_n^{i+1}}{i+1}).$$
Therefore, the restriction map is a monomorphism. To prove that the restriction map is a surjection we have to show that for a given map $\v'\in \End_{K[\der_n]}(K[x_n])$ the extension of $\v'$, which is defined as in statement 1, is a $\ggu_n$-homomorphism. The map $ \v$ is a $K$-linear map, so we have to check that $\v X_{\alpha , i } = X_{\alpha , i} \v$ for all $i=1, \ldots , n$ and $\alpha \in \N^{i-1}$.

{\em Case 1}: $i<n$. For all $\beta \in \N^{n-1}$ and $j\in \N$,
\begin{eqnarray*}
 \v X_{\alpha , i} (x^\beta \cdot x^j_n) &=& \v (x^\alpha \der_i *x^\beta x_n^j) =  \v ( \beta_ix^{\alpha +\beta - e_i} x_n^j) =   \beta_i X_{\alpha + \beta - e_i, n} \v' (\frac{x_n^{j+1}}{j+1}), \\
X_{\alpha , i} \v ( x^\beta \cdot x_n^j) &=& x^\alpha \der_i* (x^\beta \der_n * \v'(\frac{x_n^{j+1}}{j+1}))= \beta_i x^{\alpha + \beta -e_i} \der_n * \v'(\frac{x_n^{j+1}}{j+1})+x^\beta \der_n*x^\alpha\der_i*\v'(\frac{x_n^{j+1}}{j+1})\\
&=&\beta_i x^{\alpha + \beta - e_i} \der_n * \v'(\frac{x_n^{j+1}}{j+1})+0= \beta_i X_{\alpha + \beta - e_i, n}  \v'(\frac{x_n^{j+1}}{j+1}).
\end{eqnarray*}
{\em Case 2}: $i=n$. Let  $\beta \in \N^{n-1}$ and $j\in \N$. Suppose that $j\geq 1$, then
\begin{eqnarray*}
 \v X_{\alpha , n} (x^\beta  x^j_n) &=& \v (x^{\alpha +\beta} j x_n^{j-1}) =   X_{\alpha + \beta , n} \v' (x_n^j), \\
X_{\alpha , n} \v ( x^\beta  x_n^j) &=& x^\alpha \der_n* (x^\beta \der_n * \v'(\frac{x_n^{j+1}}{j+1}))= x^{\alpha +\beta}\der_n *\der_n * \v'(\frac{x_n^{j+1}}{j+1})\\
&=&X_{\alpha +\beta , n} \v' (\der_n*\frac{x_n^{j+1}}{j+1})=X_{\alpha +\beta , n}\v' (x_n^j).
\end{eqnarray*}
Suppose that $j=0$, then
\begin{eqnarray*}
 \v X_{\alpha , n} (x^\beta ) &=& \v (0) = 0, \\
X_{\alpha , n} \v ( x^\beta ) &=& x^\alpha \der_n* x^\beta \der_n * \v'(x_n)= x^{\alpha +\beta}\der_n *\der_n * \v'(x_n)\\
&=&X_{\alpha +\beta , n} \v' (\der_n*x_n)=X_{\alpha +\beta , n}\v' (1)=0 \;\; ({\rm since}\;\; \v' (1)\in K).
\end{eqnarray*}
2. Statement 2 follows from statement 1.

3. Statement 3 follows from statement 1 and Lemma \ref{a26Jan12}.

4. Statement 4 follows from statement 1. $\Box $


{\bf Monomial subspaces of the polynomial algebra $P_n$}.
Let $S$ be a subset of $\N^n$, the vector space $P_n(S):=\oplus_{\alpha \in S}Kx^\alpha$ is called a {\em monomial} subspace of the polynomial alegbra $P_n$ with {\em support} $S$. By definition, let $P_n(\emptyset ) :=0$. Clearly, $\cap_{i\in I} P_n(S_i) = P_n(\cap_{i\in I}S_i)$.


{\em Example}. For each natural number $i=1, \ldots , n$ and an ordinal number $\l \in [1,\o^n)$, the vector space $\{ p\in P_n\, | \,\frac{\der p}{\der x_i}\in P_{\l , n}\}$ is a monomial subspace of $P_n$, hence so is their intersection
$$P_{\l , n}':= \{ p\in P_n\, | \,\frac{\der p}{\der x_i}\in P_{\l , n}\;\;{\rm for}\;\; i=1, \ldots , n\}=\oplus \{ Kx^\alpha \, | \, \alpha \in \N^n,\frac{\der x^\alpha}{\der x_i}\in P_{\l , n}\;\;{\rm for}\;\; i=1, \ldots , n\}. $$
Clearly, $P_{\l , n}\subseteq P_{\l , n}'$ The next theorem describes the vector space $P_{\l , n}'$  and shows that the inclusion is always strict and $\dim_K(P_{\l , n}'/P_{\l , n})<\infty$.
\begin{theorem}\label{23Feb12}
Let $\l \in [1, \o^n)$, i.e.,  $\l = \alpha_n\o^{n-1}+\alpha_{n-1}\o^{n-2}+\cdots + \alpha_i\o^{i-1} +\cdots +\alpha_j\o^{j-1}$ with $\alpha_i\in \N$,   $\alpha_j\neq 0$ and $j\geq 1$. Then
\begin{enumerate}
\item $P_{\l , n}'=P_{\l , n}\oplus \bigoplus_{i=j}^nK\th_i $ where $ \th_i =\begin{cases}
x^\alpha =\prod_{k=1}^nx_k^{\alpha_k}& \text{if }i=j, \\
x_i\prod_{k=i}^nx_k^{\alpha_k}& \text{if }j<i\leq n.
\end{cases}$
\item $ 1\leq \dim_K(P_{\l , n}'/P_{\l , n})=n-j+1\leq n$.
\item $\l$ is not a limit ordinal (i.e.,  $\alpha_1\neq 0$)  iff $ \dim_K(P_{\l , n}'/P_{\l , n})=n$.
    \item $1\leq \dim_K(P_{\l +1,n}'/P_{\l , n}')=j\leq n-1$ and the set of elements $\{ x_ix^\alpha +P_{\l , n}'\, | \, i=1, \ldots , j\}$ is a basis for the vector space $P_{\l +1,n}'/P_{\l , n}'$.
        \item The vector spaces $\{ P_{\l , n}'\, | \, \l \in [1,\o^n )\}$ are distinct. In particular, if $\l <\mu$ then $P_{\l , n}'\varsubsetneqq P_{\mu , n}'$.
\end{enumerate}
\end{theorem}

{\it Proof}. 1. Let $R$ be the RHS of the equality in statement 1. Then $P_{\l , n}'\supseteq R$ (use Lemma \ref{a26Feb12}). In particular, the set $P_{\l , n}'\backslash  P_{\l , n}$ is a non-empty set.

{\em Case 1: $\l$ is not a limit ordinal, i.e.,  $j=1$ ($\alpha_1\neq 0$)}. In this case,
$$P_{\l , n}=\oplus \{ Kx^\alpha \, | \, x^\beta \leq x^{\alpha'} \}$$ (see (\ref{Pln1})) where $x^{\alpha'}:= x_1^{-1}x^\alpha$,  that is
$\alpha_1'=\alpha_1-1$ and $\alpha_i'=\alpha_i$ for $i=2, \ldots , n$. Let $ x^\beta \in P_{\l , n}'\backslash P_{\l , n}$. Then  $x^\beta >x^{\alpha'}$ and there exists a natural number $i$ such that $1\leq i \leq n$ and $\beta_j =\alpha_j$ for all $j>i$ and $ \beta_i = \alpha_i'+1$. Then necessarily $\beta_k=0$ for all $k<i$ (since $\frac{\der x^\beta}{\der x_k}\leq x^{\alpha'}$ for all $k<i$) and so $x^\beta =\th_i$.

{\em Case 2: $\l$ is  a limit ordinal, i.e.,  $j>1$ ($\alpha_1=0$)}.  In this case,  $$P_{\l , n} = \oplus \{ Kx^\g \, | \, x^\g <x^\alpha \}$$ (see (\ref{Pln1})).  Let $x^\beta \not\in P_{\l , n}$. Then $x^\beta \in P_{\l , n}'$ iff $\frac{\der x^\beta}{\der x_i}<x^\alpha$ for all $i=1, \ldots , n$ (recall that $P_{\l , n} = \oplus \{ Kx^\g \, | \, x^\g <x^\alpha \}$) iff

(i) $\frac{\der x^\beta}{\der x_i}\leq x^\alpha$ for all $i=1, \ldots , n$; and

(ii) $\frac{\der x^\beta}{\der x_i}\not\in K^* x^\alpha$ for all $i=1, \ldots , n$.

The condition (i) is Case 1 for $\l +1$. Therefore,
\begin{eqnarray*}
 P_{\l , n}'&=& P_{\l , n}'\cap P_{\l +1 , n}'=P_{\l , n}'\cap (P_{\l +1 , n} \oplus \bigoplus_{i=1}^nK\th_i')=P_{\l , n}'\cap (P_{\l  , n} \oplus Kx^\alpha \oplus \bigoplus_{i=1}^nK\th_i')\\
 &=&  P_{\l  , n} \oplus Kx^\alpha \oplus P_{\l , n}'\cap \bigoplus_{i=1}^nK\th_i'= P_{\l  , n} \oplus Kx^\alpha \oplus  \bigoplus_{i=1}^n P_{\l , n}'\cap K\th_i',
\end{eqnarray*}
 where
$$ \th_i':=\begin{cases}
x_i\prod_{k=i}^nx_k^{\alpha_k} & \text{if }i>j, \\
x_ix^\alpha & \text{if }i\leq j.\\
\end{cases}$$
Notice that $\th_j=x^\alpha$ and $\th_i'=\th_i\in P_{\l , n}'$ for $j<i\leq n$.
The condition (ii) excludes precisely the elements $\{ \th_i'\, | \, i\leq j\}$ that is $\oplus_{i=1}^n(P_{\l , n}'\cap K\th_i')=\oplus_{j<i\leq n}K\th_i'=\oplus_{j<i\leq n}K\th_i$. The proof of statement 1 is complete.

2. Statement 2 follows from statement 1.

3. Statement 3 follows from statement 2.

4. The ordinal number $\l +1$ is not a limit ordinal. By statement 3, $\dim_K(P_{\l +1, n}'/P_{\l +1 , n})=n$. Notice that $P_{\l , n} \subset P_{\l +1, n}\subset P_{\l +1, n}'$ and $\dim_K(P_{\l +1, n}/P_{\l , n}) =1$. By statement 1,
$$ \dim_K(P_{\l +1, n}'/P_{\l  , n}) = \dim_K(P_{\l +1, n}'/P_{\l +1 , n}) +\dim_K(P_{\l +1, n}/P_{\l  , n}) = n+1.$$
Finally, by statement 2,
$$ \dim_K(P_{\l +1, n}'/P_{\l  , n}') = \dim_K(P_{\l +1, n}'/P_{\l  , n}) -\dim_K(P_{\l , n}'/P_{\l  , n}) = n+1-(n-j+1)=j, $$
and $1\leq j=\dim_K(P_{\l +1, n}'/P_{\l  , n}') \leq n-1$. The elements $\th_i':= x_ix^\alpha$, $i=1, \ldots , j$ are the elements $\th_i$, $i=1, \ldots , j$ in statement 1 but for the ordinal $\l +1$ rather than $\l$. Clearly, $\th_i'\in P_{\l +1, n}'\backslash P_{\l , n}'$ for $i=1, \ldots , j$. Therefore, the elements $\{ \th_i'+P_{\l , n}\, | \, i=1, \ldots , j\}$  are $K$-linearly independent in the vector space $P_{\l +1, n}'/P_{\l , n}'$ as the vector spaces $P_{\l +1, n}'$ and $P_{\l , n}'$ are monomial. These elements form a basis for the vector space $P_{\l +1, n}'/P_{\l , n}'$ since $j=\dim_K(P_{\l +1, n}'/P_{\l , n}')$.

5. Statement 5 follows from statement 4.  $\Box $


For example, for all positive integers $i$ and $j$, $P_{i\o^{n-1}, n}'=P_{i\o^{n-1}, n}\oplus Kx_n^i$ and $P_{i\o^{n-1}+j, n}'=P_{i\o^{n-1}+j, n}\oplus Kx_1^jx_n^i\oplus \bigoplus_{2\leq k \leq n}K x_k x_n^i$.

{\em Example}. For each natural number $i=1, \ldots , n$, an element $\alpha \in \N^{i-1}$ and an ordinal number $\l \in [1, \o^n)$, the vector space $\{ p\in P_n\, | \, x^\alpha\frac{\der p }{\der x_i}\in P_{\l , n}\}$ is a monomial subspace of $P_n$, hence so is the intersection
$$ P_{\l ,n}'':= \{ p\in P_n \, | \; p\;\; {\rm satisfies}\;\; (\ref{Pipxn})\} = \oplus \{ Kx^\beta \, | \, \beta \in \Nn, x^\beta \;\; {\rm satisfies}\; \; (\ref{Pipxn})\}$$
where 
\begin{equation}\label{Pipxn}
P_{i-1}\frac{\der p }{\der x_i}\subseteq  P_{\l , n}\;\; {\rm for}\;\; i=1,\ldots , n.
\end{equation}
For all ordinal numbers $\l \in [1, \o^n)$,
$$P_{\l , n}\subseteq P_{\l , n}''\subseteq P_{\l , n}'.$$ The first inclusion follows from the fact that $P_{\l , n}$ is a $\ggu_n$-module. If $\l \leq \mu$ then $P_{\l , n}''\subseteq P_{\mu , n}''$. The next corollary shows that the inclusions are strict and gives a $K$-basis for every vector space $P_{\l , n}''$.
\begin{corollary}\label{a25Feb12}
Let $\l \in [1, \o^n)$, i.e.,  $\l= \alpha_m\o^{m-1}+\alpha_{m-1}\o^{m-2}+\cdots +\alpha_j\o^{j-1}$ with $\alpha_i\in \N$, $\alpha_m\neq 0$, $\alpha_j\neq 0$, $m\leq n$ and $j\geq 1$. Then
\begin{enumerate}
\item $P_{\l , n}''=P_{\l , n}\oplus Kx^\alpha =P_{\l +1, n}$ where $x^\alpha =\prod_{k=j}^mx_k^{\alpha_k}$.
\item The vector spaces $\{ P_{\l , n}''\, | \, \l \in [1, \o^n)\}$ are distinct. In particular, if $\l <\mu$  then $ P_{\l , n}''\varsubsetneqq P_{\mu , n}''$.
\item $\dim_K(P_{\l , n}''/P_{\l , n})=1$.
\item $\dim_K(P_{\l +1,n}''/P_{\l , n}'')=1$.
\end{enumerate}
\end{corollary}

{\it Proof}. 1. Statement 1 follows at once from Lemma \ref{a26Feb12}, the inclusions $P_{\l , n}\subseteq P_{\l , n}''\subseteq P_{\l , n}'$ and Theorem \ref{23Feb12}.

2-4. Statements 2-4 follow from statement 1.  $\Box $



\section{The Lie algebras $\ggu_n$ are locally finite dimensional and locally nilpotent }\label{CILAU}

The aim of this section is to prove the statement in the title of this section (Theorem \ref{22Dec11}). The key ideas are to use  the fact that the algebra $\ggu_n$ is uniserial, induction on the ordinals $\l \in [ 1, \ord (\O_n)]$ and Theorem \ref{21Dec11} that gives sufficient conditions for a Lie algebra to be  a nilpotent Lie algebra.

\begin{theorem}\label{21Dec11}
Let $J$ be an ideal of a Lie algebra $\CG$ such that the Lie factor algebra $\oCG := \CG / J$ is a finite dimensional nilpotent Lie algebra, $J$ is a nilpotent Lie algebra and every element $a\in \CG$ acts nilpotently on $J$ (that is, $(\ad \, a)^n (J)=0$ for some natural number $n=n(a)$). Then $\CG$ is a nilpotent Lie algebra.
\end{theorem}

{\it Proof}. We use induction on $d=\dim_K(\oCG )$. The case $d=0$, i.e.,  $\CG = J$, is obvious.

Let $d=1$.  This is the most important case as we will reduce the general case to this one.  Then $\CG = Ka\oplus J$ where $a\in \CG \backslash J$. Let $\d := \ad (a)\in \Inn (\CG )$ and $\der : = \ad (J) := \{ \ad (j):\CG \ra \CG  \, | \, j\in J\}$. The Lie algebra $J$ is nilpotent, that is
$$ \der^n (J) = \underbrace{\ad (J) \cdots \ad (J)}_{n \, {\rm  times}} (J)=0$$
for some natural number $n\geq 1$. For all elements $b\in J$,
$$ [ \d , \ad (b) ] = [ \ad (a) , \ad (b)]= \ad ([a,b]) \in \der, $$ hence $ \der \d \subseteq \d \der +\der$. By the assumption, the map $\d$ acts nilpotently on the ideal $J$. So, enlarging (if necessary) the number $n$ we can assume that $\d^n (J)=0$. To prove that $\CG$ is a nilpotent Lie algebra  we have to show that $(K\d +\der )^{m+1}(\CG )=0$ for some natural number $m$ ($ \CG = Ka\oplus J$ implies that $\ad (\CG ) = K\d +\der$). It suffices to show that
$$ (K\d +\der )^m (J)=0$$ for some natural number $m$  since  $(K\d +\der ) (\CG )\subseteq [ \CG , \CG ] =[Ka+J,  Ka+J]\subseteq J$. We claim that it suffices to take $m = n(n+1)$. Let $m=n(n+1)$. Then
$$(K\d +\der )^m\subseteq \sum_{s=0}^m W_s$$ where $W_s$ is the set of finite linear combinations of the elements $w_s$ where $w_s$ is a word of length $m$ in the alphabet $\{ \d , \der \}$ that contains precisely $s$ elements $\der$ and $m-s$ elements $\d$ in its product. It is sufficient to show that $w_s(J)=0$ for all $s=0,1, \ldots , m$. Notice that $\d (J)\subseteq J$ and $\d \der^i(J)\subseteq \der^i(J)$ for all $i\geq 1$. For all $s$ such that $ n \leq s\leq m$, $w_s(J)\subseteq \der^s(J)=0$ since $s\geq n$ and $\der^n(J)=0$. For all $s$ such that $0\leq s<n$, the word $w_s$ has the form $\d^{n_1}\der \d^{n_2}\der\cdots \d^{n_s}\der \d^{n_{s+1}}$ where $n_1+\cdots + n_{s+1} = m-s$. At least one of the numbers $n_i$ is $\geq n$ since otherwise we have
$$ n^2= m-n\leq m-s= n_1+\cdots + n_{s+1}\leq (s+1) (n-1) \leq n(n-1), $$
a contradiction. Then $\d^{n_i}\der \cdots \d^{n_s}\der \d^{n_{s+1}}(J)\subseteq \d^{n_i} (J) =0$ since $ n_i\geq n$ and $\d^n(J)=0$. Therefore, $w_s (J)=0$.

Let $d>1$. The Lie algebra $\oCG$ is a nilpotent Lie algebra. Fix an ideal, say $\bCJ$, of  $\oCG$  such that $\dim_K(\oCG / \bCJ )=1$. Let $\pi : \CG \ra \oCG$, $g\mapsto \overline{g} := g+J$, the canonical Lie algebra epimorphism. The ideal $\CJ := \pi^{-1}(\bCJ )$ of the Lie algebra $\CG$ has codimension 1 ($\dim_K(\CG / \CJ ) =\dim_K(\oCG / \bCJ ) =1$) and $J$ is an ideal of the Lie algebra $\CJ$. The pair $(\CJ , J)$ satisfies the condition of the theorem and $\dim_K(\CJ / J)=\dim_K(\oCG ) -1<d$. By induction on $d$, $\CJ$ is a nilpotent Lie algebra. Now, the pair $(\CG , \CJ )$ is  as in the case $d=1$ considered above. Indeed, $\dim_K(\CG / \CJ ) =1$, and the element $a\in \CG\backslash \CJ$ in the decomposition $ \CG = Ka\oplus \CJ$ acts nilpotently on $\CJ$ since $\oCG$ is a  finite dimensional nilpotent Lie algebra. In more detail, let $\d := \ad (a)$. Then $ \d^s (\CG ) \subseteq J$ for some $s$ since $\oCG$ is a finite dimensional nilpotent Lie algebra; $\d^t (J)=0$ for some $t$,  by the assumption. Hence, $\d^{s+t}(\CJ ) \subseteq \d^{s+t}(\CG ) \subseteq \d^t(J ) =0$. Therefore, $\CG$ is a nilpotent Lie algebra. $\Box $


\begin{theorem}\label{22Dec11}
The Lie algebras $\ggu_n$ are locally finite dimensional and locally nilpotent Lie algebras.
\end{theorem}

{\it Proof}. We have to show that the Lie subalgebra $\CG$ of $\ggu_n$ generated by a finite set of elements, say $a_1, \ldots , a_\nu$, of $\ggu_n$ is a finite dimensional, nilpotent Lie algebra. The case $\nu =1$ is trivial. Without loss of generality we may assume that $\nu\geq 2$, the elements $a_i$ are $K$-linearly independent and
$$ \ord (a_1)<\ord (a_2)<\cdots < \ord (a_\nu).$$
We use induction on the ordinal number $\l := \ord (a_\nu)\in [ 1, \ord (\O_n)]$. By the very definition, $\l$ is not a limit ordinal. The initial case $\l =1$ is obvious as $\CG = I_1= K\der_n$. So, let $\l$ be a non-limit ordinal such that $\l >1$, we assume that the result holds for all non-limit ordinals $\l'$ such that $\l'<\l$. Then $\l'':=\ord (a_{\nu -1})$ is a non-limit ordinal such that $\ord (a_i) \leq \l''<\l$ for all $i=1, \ldots , \nu -1$. By induction, the Lie subalgebra, say $V$, of $\ggu_n$ generated by the elements $a_1, \ldots , a_{\nu -1}$ is a finite dimensional, nilpotent Lie algebra. The inner derivation $\d : = \ad (a_\nu)$ of the Lie algebra $\ggu_n$ is a locally nilpotent derivation (Theorem \ref{a8Dec11}.(5)) and $\dim_K(V)<\infty$, hence $\d^{s+1}(V)=0$ for some natural number $s$. The vector space $U:= V+\d (V)+\cdots + \d^s (V)$ is a finite dimensional, $\d$-invariant subspace of the ideal $I_{\l''}$ (see (\ref{Jlid})).  By induction, the Lie subalgebra, say $J$, of $\ggu_n$ generated by $U$ is a finite dimensional, nilpotent Lie algebra. Then $\d^t(J)=0$ for some natural number $t\geq 1$. Clearly, $\d (J)\subseteq J$ since $\d (U)\subseteq U$ and $\d$ is a derivation. We see that $\CG = Ka_\nu +J= Ka_\nu \oplus J$ (clearly, $\CG \supseteq Ka_\nu +J$; on the other hand, $Ka_\nu +J$ is a Lie subalgebra of $\ggu_n$ that contains the elements $a_1, \ldots , a_\nu $, and so  $\CG \subseteq Ka_\nu +J$). Claim: {\em The pair $(\CG , J)$ satisfies the assumptions of Theorem \ref{21Dec11}}, hence $\CG$ is a nilpotent Lie algebra (by Theorem \ref{21Dec11})), and $\dim_K(\CG ) = 1+\dim_K(J)<\infty$.  To prove the claim  it suffices to show that $(\d + \der )^m(J)=0$ for all $m\gg 0$ where $\der = \ad (J)$. Fix a number $n$ such that $\d^n (J) =0$ and $\der^n (J)=0$. By a similar reason as in the proof of Theorem \ref{21Dec11} it suffices to take $m=n(n+1)$:
$$ (\d +\der )^m(J)\subseteq \sum_{s=0}^{n-1} \sum_{n_1+\cdots n_{s+1}=m-s}\d^{n_1}\der \d^{n_2}\der\cdots \d^{n_s}\der \d^{n_{s+1}}(J)=0. \;\;\; \Box $$



\section{The  isomorphism problem for the factor algebras of the Lie  algebras $\ggu_n$}\label{ISOPRUN}
We know already that the Lie algebras $\ggu_n$ and $\ggu_m$ are not isomorphic for $n\neq m$.
The aim of this section is to answer the question (Theorem \ref{24Dec11} and Corollary \ref{a24Dec11}): {\em Let $I$ and $J$ be ideals of the Lie algebras $\ggu_n$ and $\ggu_m$ respectively. When the Lie algebras $\ggu_n/I$ and $\ggu_m/J$ are isomorphic?} First, we consider the case when $n=m$ (Theorem \ref{24Dec11}) and then the general case (Corollary \ref{a24Dec11}) will be deduced from the this one. Let $\o^0:=1$.

\begin{theorem}\label{24Dec11}

\begin{enumerate}
\item Let $I$ be an  ideal of the Lie algebra $\ggu_n$, that is $I=I_\l$ for some $\l \in [1, \ord (\O_n)]\cup \{ 0\}$ (Theorem \ref{10Dec11}) where $I_0:=\{ 0\}$. Then the Lie algebras $\ggu_n$ and $\ggu_n/I_\l$ are isomorphic iff $\l = i\o^{n-2}$ where $i\in \N $.
\item Let $I$ and $J$ be  ideals of the Lie algebra $\ggu_n$, that is $I=I_\l$ and $J=I_\mu$ for some elements $\l , \mu \in [ 1, \ord (\O_n)]\cup \{ 0\}$ (Theorem \ref{10Dec11}.(1)). Then the Lie algebras $\ggu_n / I_\l$ and $\ggu/I_\mu$ are isomorphic iff
    \begin{enumerate}
\item $\l = i\o^{n-2}+\nu$ and $\mu = j\o^{n-2}+\nu$ where $i,j\in \N$ and $\nu \in [1,\o^{n-2})\cup \{ 0\}$; or
\item $\l = \o^{n-1}+\o^{n-2}+\cdots +\o^s+i\o^{s-2}+\nu $ and $\mu = \o^{n-1}+\o^{n-2}+\cdots +\o^s+j\o^{s-2}+\nu $ where $2\leq s\leq n-1$; $i,j\in \N$ and $\nu \in [1,\o^{s-2})\cup\{ 0\}$; or
\item $\l = \mu = \o^{n-1}+\o^{n-2}+\cdots +\o+\varepsilon $ where $\varepsilon =0,1$.
\end{enumerate}
\end{enumerate}
\end{theorem}

For natural numbers $n$ and $m$ such that $2\leq n<m$, there is a natural Lie algebra isomorphism
\begin{equation}\label{uunm}
\ggu_n \simeq \ggu_m / I_{\nu_{mn}}, \;\; X_{\alpha , i}\mapsto X_{\alpha , i}+I_{\nu_{mn}},\;\; {\rm  where} \;\; X_{\alpha , i}\in \CB_n, \;\;\nu_{mn}:= \o^{m-1}+\o^{m-2}+\cdots + \o^n.
\end{equation}
In more detail, $I_{\nu_{mn}}=\ggu_{m,n+1}=\oplus_{j=n+1}^mP_{j-1}\der_j$.
\begin{corollary}\label{a24Dec11}
Let $n$ and $m$ be natural numbers such that $2\leq n <m$, $I$ and $J$ be ideals of  the Lie algebras $\ggu_n$ and $\ggu_m$ respectively. Then the Lie algebras $\ggu_n / I$ and $\ggu_m / J$ are isomorphic iff
\begin{eqnarray*}
 (I,J)&\in &  \{ (I_\l , I_\mu ) \, | \, \l =i\o^{n-2}+\nu, \mu =\o^{m-1}+\o^{m-2}+\cdots +\o^n+j\o^{n-2}+\nu \\
 & &  {\rm where}\;\; \nu\in [ 1, \o^{n-2})\cup \{ 0\} \;\; {\rm and}\;\; i,j\in \N\}\\
 &\cup & \bigcup_{s=2}^{n-1}\{ (I_\l , I_\mu ) \, | \, \l = \o^{n-1}+\o^{n-2}+\cdots +\o^s+ i\o^{s-2}+\nu, \mu =\o^{m-1}+\o^{m-2}+\cdots +\o^s+\\
  & & +j\o^{s-2}+\nu\;\;  {\rm where}\;\; \nu\in [ 1, \o^{s-2})\cup \{ 0\}\;\; {\rm and}\;\; i,j\in \N\}\\
  &\cup& \{ (I_\l , I_\mu )\, | \, \l = \o^{n-1}+\o^{n-2}+\cdots +\o+ \varepsilon , \mu =\o^{m-1}+\o^{m-2}+\cdots +\o +\varepsilon  \\
   & &  {\rm where}\;\;  \varepsilon =0,1\}.
  \end{eqnarray*}
\end{corollary}


{\bf The Lie algebra epimorphism $f_n:\ggu_n\ra\ggu_n$ with nonzero kernel}. Recall that $\ggu_n=\ggu_{n-1}\oplus P_{n-1}\der_n$ where $\ggu_{n-1}$ is a Lie subalgebra of $\ggu_n$ and $P_{n-1}\der_n$ is an abelian ideal of the Lie algebra $\ggu_n$. Earlier, we introduced the $K$-basis $\CB_n=\{ X_{\alpha , i}\}$ for the Lie algebra $\ggu_n$ (see (\ref{Xai})). Let us define the $K$-linear map $f_n:\ggu_n\ra \ggu_n$,  by one of the two  equivalent ways:
\begin{equation}\label{fnun}
  f_n(u):=\begin{cases}
u& \text{if }u\in \ggu_{n-1},\\
[\der_{n-1},u]& \text{if }u\in P_{n-1}\der_n, \\
\end{cases}
\;\;\; f_n(X_{\alpha , i}):=\begin{cases}
X_{\alpha , i}& \text{if }\alpha\in \N^{i-1}, \; i\neq n,\\
\alpha_{n-1}X_{\alpha -e_{n-1}, n}& \text{if }\alpha\in \N^{n-1}, \; i= n. \\
\end{cases}
\end{equation}
\begin{lemma}\label{ap10Dec11}
For every natural number $i\geq 1$, the $K$-lineal map $f_n^i:\ggu_n\ra \ggu_n$ is a Lie algebra epimorphism and $\ker (f_n^i) = \sum_{j=0}^{i-1}P_{n-2}x_{n-1}^j\der_n=I_{i\o^{n-2}}$.
\end{lemma}

{\it Proof}. By the very definition of the map $f_n$,
$$\ker (f_n^i) = \ker_{P_{n-1}\der_n}(\ad (\der_{n-1})^i)= \sum_{j=0}^{i-1}P_{n-2}x_{n-1}^j\der_n=I_{i\o^{n-2}}.$$
To finish the proof of the lemma it suffices to show that the $K$-linear map $f_n$ is a Lie algebra homomorphism, that is $f_n([u,v])=[f_n(u), f_n(v)]$ for all elements $u$ and $v$ of the basis $\CB_n$ of the Lie algebra $\ggu_n$ (see (\ref{Xai})). Since $\ggu_n = \ggu_{n-1}\oplus P_{n-1}\der_n$, we have $\CB_n   = \CB_{n-1}\cup\{ X_{\alpha , n}\}_{\alpha \in \N^{n-1}}$. The equality above is obvious if either $u,v\in \CB_{n-1}$ (since $f_n(a) = a$ for all $a\in \ggu_{n-1}$) or $u,v\in \{ X_{\alpha , n}\}_{\alpha \in \N^{n-1}}$ (since $f_n(P_{n-1}\der_n) \subseteq P_{n-1}\der_n$ and $[P_{n-1}\der_n,P_{n-1}\der_n]=0$). In the remaining case where $u\in \CB_{n-1}$ and $v=X_{\alpha , n}$ for some $\alpha \in \N^{n-1}$, the equality follows from the following three facts: $\der_{n-1}$ is a central element of the Lie algebra $\ggu_{n-1}$ (Proposition \ref{a8Dec11}.(6)), $f_n(P_{n-1}\der_n)\subseteq P_{n-1}\der_n$ and $f_n|_{P_{n-1}\der_n }= \ad (\der_{n-1})$. Indeed, applying the inner derivation $\d := \ad (\der_{n-1})$  of the Lie algebra $\ggu_n$ to the equality $[u,X_{\alpha , n}]=w$ (where $w\in P_{n-1}\der_n$) we obtain the required equality:
\begin{eqnarray*}
 f_n(w)&=& \d([u,X_{\alpha , n}])=[\d (u),X_{\alpha , n}]+[u,\d (X_{\alpha , n})]=[0,X_{\alpha , n}]+[f_n(u),f_n(X_{\alpha , n})]\\
 & =& [f_n(u),f_n(X_{\alpha , n})].\;\; \Box
\end{eqnarray*}


\begin{corollary}\label{ya24Jan12}
Let $n\geq 3$. The ideal $P_{n-1}\der_n$ is the least ideal $I$ of the Lie algebra $\ggu_n$ such that the Lie factor algebra $\ggu_n / I$ is isomorphic to the Lie algebra $\ggu_{n-1}$.
\end{corollary}

{\it Proof}. It is obvious that $\ggu_n/P_{n-1}\der_n\simeq \ggu_{n-1}$ and $P_{n-1}\der_n = I_{\o^{n-1}}$. Suppose that $I\neq P_{n-1}\der_n$, we seek a contradiction. Then $I=I_\l$ for some  $\l < \o^{n-1}$. Fix a natural number $i\geq 1$ such that $\l <i\o^{n-2}$. Then $I_\l \subset I_{i\o^{n-2}}$.  There is a natural  epimorphism of Lie algebras $\ggu_n/I\ra \ggu_n/ I_{i\o^{n-2}}$, and so $ \udim (\ggu_n/ I_{i\o^{n-2}})\leq \udim (\ggu_n/I)=\udim (\ggu_{n-1})$. By Lemma \ref{ap10Dec11}, $\ggu_n/I_{i\o^{n-2}}\simeq \ggu_n$, hence
$$\udim (\ggu_{n-1}) \geq \udim ( \ggu_n / I_{i\o^{n-2}})=\udim (\ggu_n) >\udim (\ggu_{n-1}),$$
a contradiction.  $\Box $

$\noindent $

{\bf Proof of Theorem \ref{24Dec11}}. 1. $(\Leftarrow )$ By Lemma \ref{ap10Dec11}, for every natural number $i\geq 1$, the map $f_n^i:\ggu_n\ra \ggu_n$ is a Lie algebra epimorphism with kernel $I_{i\o^{n-2}}$. Therefore, $\ggu_n\simeq \ggu_n/I_{i\o^{n-2}}$.

$(\Rightarrow )$ This implication follows from statement 2.

2. To prove statement 2 we use induction on $n\geq 2$. The initial step $n=2$ is a direct consequence  of Lemma \ref{ap10Dec11} and the classification of the ideals of the Lie algebra $\ggu_2$. By Theorem \ref{10Dec11}.(1), the proper ideals of the Lie algebra $\ggu_2$ are $I_n$ ($n\geq 1$) and $I_\o = P_1\der_2$. By Lemma \ref{ap10Dec11}, $\ggu_2/I_n\simeq \ggu_2$ for all $n\geq 1$; $\dim_K(\ggu_2/I_\o ) = \dim_K(K\der_1) = 1<\dim_K(\ggu_2)=\infty$; and statement 2 follows.

Let $n>2$ and we assume that the result holds for all $n'<n$. Recall that $\ggu_{n-1}\simeq \ggu_n/I_{\o^{n-1}}$, $I_{\o^{n-1}}=P_{n-1}\der_n$ and $\udim (\ggu_{n-1}) < \udim (\ggu_n)$. The Lie algebra $\ggu_n$ is a uniserial Artinian Lie algebra, hence so is each of its Lie factor algebras. So, if $ \ggu_n / I_\l \simeq \ggu_n / I_\mu$ then $\udim (\ggu_n / I_\l)= \udim ( \ggu_n / I_\mu )$.

{\em Step 1: $\udim (\ggu_n / I_\l ) =\udim (\ggu_n)$ for all} $\l <\o^{n-1}$.

Indeed, $\l < \o^{n-1}$ implies $\l \leq i\o^{n-2}$ for some $i\geq 1$. Then $\ggu_n \simeq \ggu_n / I_{i\o^{n-2}}$ (Lemma \ref{ap10Dec11}), and so
$$ \udim (\ggu_n ) \geq \udim (\ggu_n / I_\l ) \geq \udim ( \ggu_n / I_{i\o^{n-2}})=\udim (\ggu_n), $$
 hence $\udim (\ggu_n ) = \udim (\ggu_n / I_\l )$.

 Suppose that $\ggu_n/ I_\l \simeq \ggu_n / I_\mu$ for some $\l $ and $\mu$.

 {\em Step 2: It is sufficient to consider the case when $ \l, \mu <\o^{n-1}$}.

 If $\l , \mu \geq \o^{n-1}$ then statement 2 follows by induction on $n$ as $I_{\o^{n-1}}\subseteq I_\l$, $I_{\o^{n-1}}\subseteq I_\mu$ and $\ggu_n / I_{\o^{n-1}}\simeq \ggu_{n-1}$. Without loss of generality we may assume that $\l \leq \mu$. The case when $\l < \o^{n-1} \leq \mu$ is not possible, by Step 1,  since
 $$ \udim (\ggu_n / I_\l )=\udim (\ggu_n)  >\udim (\ggu_{n-1}) = \udim (\ggu_n/ I_{\o^{n-1}})\geq \udim (\ggu_n / I_\mu )=\udim (\ggu_n / I_\l ),$$
  a contradiction.  This finishes the proof of Step 2.

{\em Step 3: In view of Lemma \ref{ap10Dec11}, we may assume that $\l , \mu <\o^{n-2}$}.

The idea of the proof of the theorem is, for each Lie factor algebra $\ggu_n/I_\l$ where $\l <\o^{n-2}$, to introduce an isomorphism invariant that is {\em distinct} for all ordinal numbers $\l <\o^{n-2}$. The invariant is the uniserial dimension of certain ideals of the Lie algebra $\ggu_n/I_\l$.  For each ordinal number $\nu \in [ 1, \udim (\ggu_n / I_\l )]$, let $I_\nu'$ be the unique ideal of the uniserial Artinian Lie algebra $\ggu_n/I_\l$ with uniserial dimension $\nu$. Every ideal of uniserial Artinian Lie algebra is characteristic ideal. This is the  crucial fact in the arguments below.   Clearly,
\begin{eqnarray*}
 I_{\o^{n-2}}'&\equiv &I_{\o^{n-2}}\equiv P_{n-2}\der_n \mod I_\l ,  \\
 I_{\o^{n-1}+\o^{n-2}}'&\equiv &I_{\o^{n-1}+\o^{n-2}}\equiv P_{n-2}\der_{n-1}+P_{n-1}\der_n \mod I_\l .
\end{eqnarray*}
Let $ \ga := I_{\o^{n-1}+\o^{n-2}}'$ and $\th := x_{n-1}\der_n +I_\l \in \ggu_n / I_\l$. Notice that $I_{\o^{n-2}+1}'=K\th \oplus I_{\o^{n-2}}'$. For an element $v\in \ggu_n / I_\l$, let $\Cen_\ga (v) := \{ a\in \ga \, | \, [a,v]=0\}$.

{\em Step 4: For all $\s \in \Aut_K(\ggu_n/ I_\l )$ and $ u\in S:= I_{\o^{n-2}+1}'\backslash I_{\o^{n-2}}'$, $\Cen_\ga (\s (u))=\Cen_\ga (\th ) = I_{\o^{n-1}+\l}'$.}

The second equality follows from  two facts : $[x_{n-1}\der_n, I_{\o^{n-1}}]=[x_{n-1}\der_n, P_{n-1}\der_n]=0$ and, for all $\alpha\in \N^{n-2}$, $[x^\alpha \der_{n-1},  x_{n-1}\der_n]=x^\alpha \der_n$.

The first equality follows from the facts that $\s (u) = \xi \th +v+I_\l$ for some elements $\xi \in K^*$,  $v\in P_{n-2}\der_n$ and $$[v,I_{\o^{n-1}+\o^{n-2}}]\subseteq [P_{n-2}\der_n, P_{n-2}\der_{n-1}+P_{n-1}\der_n]=0.$$ Alternatively, using the equality $ \Cen_\ga ( \th ) = I_{\o^{n-1}+\l}'$ and the fact that $I_{\o^{n-1}+\l}'$ and $\ga$ are characteristic ideals of the Lie algebra $\ggu_n/I_\l$ we see that
$$ I_{\o^{n-1}+\l}'=\s (I_{\o^{n-1}+\l}')=\s (\Cen_\ga (\th ))= \Cen_{\s (\ga )} (\s (\th ) )=\Cen_\ga (\s (\th )).$$
Since $\s (\ga ) = \ga$ and $\s (S)=S$ for all automorphisms $\s \in \Aut_K(\ggu_n/I_\l ) $, Step 4 means that the ordinal number
$$\o^{n-1}+\l = \udim (I_{\o^{n-1}+\l}')= \udim (\Cen_\ga (\s (u)))$$ is an isomorphism invariant for the algebra $\ggu_n/I_\l$ (where $\udim (I_{\o^{n-1}+\l})$ is the uniserial dimensions of the $\ggu_n$-{\em module} $I_{\o^{n-1}+\l}$). If $\ggu_n/I_\l \simeq \ggu_n/I_\mu$ for some ordinals $\l ,\mu <\o^{n-2}$ then  $\o^{n-1}+\l = \o^{n-1}+\mu$, hence $\l = \mu$. $\Box $


\begin{corollary}\label{xx24Dec11}
Let $I_\l$ be an ideal of the Lie algebra $\ggu_n$ where $\l \in [1, \ord (\O_n)]\cup \{ 0 \}$ where $I_0:=0$. Then
\begin{enumerate}
\item $\udim (\ggu_n / I_\l) = \udim (\ggu_n)$ for all $\l \in [1, \o^{n-1})\cup \{ 0\}$,
\item $\udim (\ggu_n / I_\l) = \udim (\ggu_s)$ for all $\l \in [\o^{n-1}+\o^{n-2}+\cdots +\o^s, \o^{n-1}+\o^{n-2}+\cdots +\o^{s-1})$ and $s$ such that $2\leq s \leq n-1$.
\item $\udim (\ggu_n / I_\l) =1-\epsilon $ for all $\l= \o^{n-1}+\o^{n-2}+\cdots +\o+\epsilon$ where $\epsilon = 0,1$.
\end{enumerate}
\end{corollary}

{\it Proof}. 1. Statement 1 is Step 1 in the proof of Theorem \ref{24Dec11}.

2. Statement 2 follows from statement 1 and (\ref{uunm}).

3. Trivial. $\Box $

$\noindent $

{\bf Proof of Corollary \ref{a24Dec11}}. Notice that the uniserial dimensions $\udim (\ggu_i)$, $i\geq 2$, are distinct, and if $\ggu_n/I_\l\simeq \ggu_m/ I_\mu$ then $\udim (\ggu_n/I_\l )=\udim ( \ggu_m/ I_\mu)$.  The Lie algebra isomorphism $\ggu_n \simeq \ggu_m / I_{\nu_{mn}}$ (see (\ref{uunm})) induces the bijection from the set $\CJ (\ggu_n)$ of all the nonzero ideals of the Lie algebra $\ggu_n$  to the set $\CJ (\ggu_m , I_{\nu_{mn}})$ of all the ideals of the Lie algebra $\ggu_m$ that {\em properly} contain the ideal $I_{\nu_{mn}}$:
\begin{equation}\label{CJmn}
\CJ (\ggu_n)\ra \CJ (\ggu_m , I_{\nu_{mn}}), \;\; I_\l \mapsto I_{\nu_{mn}+\l }.
\end{equation}
Now, the corollary follows from Theorem \ref{24Dec11}, Corollary \ref{xx24Dec11} and (\ref{CJmn}).  $\Box $


\section{The  Lie algebra $\ggu_\infty$}\label{LIUIN}

In this section, the Lie algebra $\ggu_\infty$ is studied in detail. Many properties of the Lie algebra $\ggu_\infty$ are similar to those of the Lie algebras $\ggu_n$ ($n\geq 2$) but there are several differences. For example, the Lie algebra $\ggu_\infty$ is not solvable, not Artinian  but almost Artinian, $\udim (\ggu_\infty ) = \o^\o$. A classification of all the ideals of the Lie algebra $\ggu_\infty$ is obtained (Theorem \ref{26Dec11}), all of them are characteristic ideals (Corollary \ref{a26Dec11}.(2)). An isomorphism criterion is given for the Lie factor algebras of $\ggu_\infty$ (Corollary \ref{ca24Dec11}).

Let $P_\infty := \cup_{n\geq 1} P_n = K[x_1, x_2, \ldots , ]$, the polynomial algebra in countably many variables, and $A_\infty := \cup_{n\geq 1} A_n = K\langle x_1, x_2, \ldots , \der_1, \der_2, \ldots \rangle$, the infinite Weyl algebra. The Lie algebra $\ggu_\infty$ is a Lie subalgebra  of the Lie algebra $(A_\infty , [ \cdot , \cdot ] )$. The polynomial algebra $P_\infty$ is an $A_\infty$-module. In particular, the polynomial algebra $P_\infty$ is a $\ggu_\infty$-module. The Lie algebra
$$u_\infty :=\bigoplus_{i\geq 1} P_{i-1}\der_i$$
is the direct   sum  of abelian (infinite dimensional for $i\geq 2$)  Lie subalgebras $P_{i-1}\der_i$. For each natural number $i\geq 1$, 
$$\ggu_{\infty , i}:= \oplus_{j\geq i} P_{j-1}\der_j$$ is an ideal of the Lie algebra $\ggu_\infty$, by (\ref{PdiPdj}). Clearly, $\ggu_{i,i}\subset \ggu_{i+1, i} \subset  \cdots \subset \ggu_{\infty , i}=\cup_{n\geq i}\ggu_{n,i}$ for all $i\geq 2$. There is the strictly descending chain of ideals of the Lie algebra $\ggu_\infty$,
\begin{equation}\label{infseru}
\ggu_{\infty , 1}= \ggu_{\infty} \supset \ggu_{\infty , 2} \supset \cdots \supset \ggu_{\infty , n} \supset \cdots \supset \bigcap_{i\geq 1} \ggu_{\infty , i}=0
\end{equation}
and $\ggu_\infty / \ggu_{\infty , n+1}\simeq \ggu_n$ for all $n\geq 2$.

\begin{proposition}\label{25Dec11}

\begin{enumerate}
\item The Lie algebra $\ggu_\infty$ is not a solvable Lie algebra.
\item  The Lie algebra $\ggu_\infty$ is a locally nilpotent  and locally finite dimensional Lie algebra.
\item Each element $u\in \ggu_\infty$ acts locally nilpotently on the $\ggu_\infty$-module $P_\infty$.
\item The chain of nonzero ideals in (\ref{infseru}) is the derived  series for the Lie algebra $\ggu_\infty$, that is $(\ggu_\infty)_{(i)}=\ggu_{\infty , i+1}$ for all $i\geq 0$.
\item The upper central series for the Lie algebra $\ggu_\infty$ stabilizers at the first step, that is $(\ggu_{\infty})^{(0)}=\ggu_{\infty}$ and $(\ggu_{\infty})^{(i)}=\ggu_{\infty , 2}$ for all $i\geq 1$.
\item All the inner derivations of the Lie algebra $\ggu_\infty$ are locally nilpotent derivations.
\item The centre $Z(\ggu_\infty )$ of the Lie algebra $\ggu_\infty$ is $0$. In particular, $\cdim (\ggu_\infty )=0$.
\item  The Lie algebras $\ggu_\infty$ and $\ggu_n$ where $n\geq 2$ are not isomorphic.
    \item The Lie algebra $\ggu_\infty$ contains a copy of every nilpotent finite dimensional Lie algebra.
        \item Let $u\in \ggu_\infty$. Then the inner derivation $\ad (u)$ is a nilpotent derivation of the Lie algebra $\ggu_\infty$     iff $a=0$.
\end{enumerate}
\end{proposition}

{\it Proof}. 1. Statement 1 follows from statement 4.

2. Statement 2 follows from Theorem \ref{22Dec11} and the fact that $\ggu_\infty = \cup_{n\geq 2}\ggu_n$.

3. Statement 3 follows from Proposition \ref{a8Dec11}.(4)  and the facts that $P_\infty = \cup_{n\geq 1} P_n$ and $\ggu_\infty =\cup_{n\geq 2}\ggu_n$.

4, 5. Statements 4 and 5  follow from (\ref{PdiPdj}) and the decomposition $\ggu_\infty = \oplus_{i\geq 1}P_{i-1}\der_i$.

6.  Statement 6 follows from Proposition \ref{a8Dec11}.(5) and the facts that  $\ggu_\infty =\cup_{n\geq 2}\ggu_n$ and $\ggu_2\subseteq \ggu_3\subseteq \cdots $.

7. If $z\in Z(\ggu_\infty )$ then $z\in \ggu_n$ for some $n$ and so $z\in Z(\ggu_n)=K\der_n$ (Proposition \ref{a8Dec11}.(6)). Since $[\der_n, x_n\der_{n+1}]=\der_{n+1}$, we must have $z=0$.

8. The Lie algebra $\ggu_n$ is solvable (Proposition \ref{a8Dec11}.(1)) but the Lie algebra $\ggu_\infty$ is not (statement 1). Therefore, $\ggu_\infty \not\simeq \ggu_n$ for all $n\geq 2$. By Proposition \ref{a8Dec11}.(7), the Lie algebras $\ggu_n$, $n\geq 2$, are pairwise non-isomorphic.

9. Any nilpotent finite dimensional Lie algebra is a subalgebra of the Lie algebra $\UT_n(K)$ for some $n\geq 2$. Now, the result follows from the inclusions $\UT_n(K)\subseteq \ggu_n\subseteq \ggu_\infty$.

10. It suffices to show that if $a\neq 0$ then the derivation $\d$ is not nilpotent. Suppose this is not the case for some $a$, we seek a contradiction. We can write the element $a$ as the sum $p_n\der_n+p_{n+1}\der_{n+1}+\cdots $ where $p_i\in P_{i-1}$ for all $i\geq n$ and $p_n\neq 0$. In view of the Lie algebra isomorphism $\ggu_\infty / \ggu_{\infty , n+2}\simeq \ggu_{n+1}$, the induced inner derivation $\d'= \ad (p_n\der_n+p_{n+1}\der_{n+1})$ by $\d$ is a nilpotent derivation of the Lie algebra $\ggu_{n+1}$. By Lemma \ref{b19Jan12}, $p_n=0$, a contradiction.  $\Box $

$\noindent $

The next theorem gives a list of all the ideals of the Lie algebra $\ggu_\infty $.

\begin{theorem}\label{26Dec11}
The set $\CJ (\ggu_\infty )$ of all the nonzero ideals of the Lie algebra $\ggu_\infty$ is equal to the set $\{ \ggu_\infty , \ggu_{\infty , n}, I_{\l , n}:= I_\l (n) +\ggu_{\infty , n+1}\, | \, n\geq 2, \l \in [1,\o^{n-1})\}$ where $I_\l ( n)$ is the ideal $I_\l$ of the Lie algebra $\ggu_n$ as defined in (\ref{Jlid}), that is $I_\l (n) = \oplus_{(\alpha , n)\leq \l} KX_{\alpha , n}$. In particular, every nonzero ideal $I$ of the Lie algebra $\ggu_\infty$ contains an ideal $\ggu_{\infty , n'+1}$ for some $n'\geq 2$, and, for $I=I_{\l , n}$, $n=\min \{ n'\geq 2\, | \, \ggu_{\infty , n'+1}\subseteq I\}$.
\end{theorem}

{\it Proof}. Let $I$ be a nonzero ideal of the Lie algebra $\ggu_\infty$ and $a$ be a nonzero element of $I$. Then $$ a=a_l\der_l+a_{l+1}\der_{l+1}+\cdots  + a_m\der_m=a_l\der_l+\cdots$$ for some elements $l\geq 1$, $a_i\in P_{i-1}$ and $a_l\neq 0$. Let $d$ be the (total) degree of the polynomial $a_l=\sum_{\alpha \in \N^{l-1}}\l_\alpha x^\alpha\in P_{l-1}$ where $\l_\alpha \in K$. Fix $\alpha = (\alpha_i) \in \N^{l-1}$ such that $|\alpha | := \alpha_1+\cdots +\alpha_{l-1}=d$. Applying $\ad (\der )^\alpha := \prod_{i=1}^{l-1}\ad (\der_i)^{\alpha_i}$ to the element $a$ yields the element of the ideal $I$ of the type $\alpha !\l_\alpha \der_l+\cdots $. So, without loss of generality we may assume from the very beginning that $a_l=1$, that is $a=\der_l+a_{l+1}\der_{l+1}+\cdots + a_m\der_m$. Then
$$ I\ni [ a,x_l\der_{m+1}]= [ \der_l , x_l\der_{m+1}]=\der_{m+1}.$$Hence, for all $s>m+1$, $I\supseteq [ \der_{m+1}, P_{s-1}\der_s]=[\der_{m+1}, P_{s-1}]\der_s= P_{s-1} \der_s$. This means that $\ggu_{\infty , m+2} \subseteq I$. Consider the Lie algebra epimorphism (where $n=\min \{ n'\geq 2\, | \, \ggu_{\infty , n'+1}\subseteq I\}$)
\begin{equation}\label{uinun}
\pi_n : \ggu_\infty \ra \ggu_\infty / \ggu_{\infty , n+1}\simeq \ggu_n, \;\; u\mapsto \bu := u+\ggu_{\infty , n+1}.
\end{equation}
The image $\pi_n (I)$ of the ideal $I$ is an ideal of the algebra $\ggu_n$ such that $\pi_n(I)\subset \ggu_{n,n}=P_{n-1}\der_n$, by the very definition of the number $n$. Without loss of generality we may assume that $I\neq \ggu_\infty, \ggu_{\infty, 2}, \ldots $ Then $I=I_{\l , n}$ for some $\l \in [1, \o^{n-1})$, by Theorem \ref{10Dec11}.(1). It follows from the two obvious facts:

(i) $I_{\l , n}\supseteq I_{\mu , m}$ iff either $n<m$ or, otherwise, $n=m$ and $\l \geq \mu$;

(ii) $\ggu_{\infty , n-1}\supset I_{\l , n}\supset \ggu_{\infty , n}$ for all $n\geq 2$;

that the Lie algebra $\ggu_\infty$ is a uniserial Lie algebra (that is, for any two distinct ideals $I$ and $J$ of $\ggu_\infty$ either $I\subset J$ or $I\supset J$). So, the chain
\begin{equation}\label{uinfid}
\ggu_\infty = \ggu_{\infty , 1}\supset \cdots \supset I_{\l , 2}\supset \cdots \supset I_{1,2}\supset \ggu_{\infty , 2}\supset \cdots \supset \ggu_{\infty, n}\supset \cdots \supset I_{\mu , n}\supset \cdots \supset I_{1, n}\supset \ggu_{\infty , n+1}\supset \cdots
\end{equation}
contains all the nonzero ideals of the Lie algebra $\ggu_\infty$. $\Box $


In combination with Theorem \ref{24Dec11} and Corollary \ref{a24Dec11}, the next corollary gives an isomorphism  criterion for the Lie  factor algebras
of the Lie algebra $\ggu_\infty$.

\begin{corollary}\label{ca24Dec11}

\begin{enumerate}
\item Let $I$ be an ideal of the Lie algebra $\ggu_\infty$. Then $\ggu_\infty / I \simeq \ggu_\infty$ iff $I=0$.
\item Let $I$ and $J$  be nonzero  ideals of the Lie algebra $\ggu_\infty$ and $n=\min \{ n'\geq 2\, | \, \ggu_{\infty , n'+1}\subseteq I, \ggu_{\infty , n'+1}\subseteq  J\}$ ($n<\infty$, see Theorem \ref{26Dec11}). Then $\ggu_\infty / I \simeq \ggu_\infty/J$ iff $\ggu_n/ I' \simeq \ggu_n/J'$ where $I':= I/\ggu_{\infty , n+1}$, $J':= J/\ggu_{\infty , n+1}\subseteq \ggu_n = \ggu_\infty / \ggu_{\infty , n+1}$.
\end{enumerate}
\end{corollary}
In contrast to the Lie algebras $\ggu_n$, $n\geq 2$, no proper Lie factor algebras of $\ggu_\infty$ is isomorphic to $\ggu_\infty$.

A Lie algebra $\CG$ is called an {\em almost Artinian} Lie algebra if all the {\em proper} factor algebras are Artinian Lie algebras (i.e.,  for every nonzero ideal $I$ of the Lie algebra $\CG$, the factor algebra $\CG / I$ is an Artinian Lie algebra). \begin{corollary}\label{a26Dec11}

\begin{enumerate}
\item The Lie algebra $\ggu_\infty$ is a uniserial, neither Artinian nor Noetherian, almost Artinian  Lie algebra, and its uniserial dimension is equal to $\udim (\ggu_\infty ) = \o^{\o}$.
\item All the ideals of the Lie algebra $\ggu_\infty$ are characteristic ideals.
\end{enumerate}
\end{corollary}

{\it Proof}. 1. We have just seen that the Lie algebra $\ggu_\infty$ is uniserial (see (\ref{uinfid})). By Theorem \ref{26Dec11}, the Lie algebra $\ggu_\infty$ is neither Artinian nor Noetherian but it is an almost Artinian Lie algebra (since all the Lie algebras $\ggu_n$ are Artinian, Theorem \ref{10Dec11}.(2)). By Theorem \ref{26Dec11}, every nonzero ideal of the Lie algebra $\ggu_\infty$ contains the ideal $\ggu_{\infty , n+1}$ for some $n$, and the factor algebra $\ggu_\infty/\ggu_{\infty , n+1}\simeq \ggu_n$ is a uniserial Artinian Lie algebra, hence
$$ \udim (\ggu_\infty ) = \max \{ \udim (\ggu_\infty / \ggu_{\infty, n+1}) = \udim  (\ggu_n) = \o^{n-1}+\o^{n-2}+\cdots + \o +1\, | \, n\geq 2\} = \o^{\o},$$
by Theorem \ref{10Dec11}.(2).

2. The ideals $\ggu_{\infty, n}$, $n\geq 1$ are characteristic ideals  as they form the derived  series of the Lie algebra $\ggu_\infty$ (Proposition \ref{25Dec11}.(4)). By Theorem \ref{26Dec11}, it  remains to show that every ideal $I_{\l , n}$ is a characteristic ideal. Let $\s$ be an automorphism of the Lie algebra $\ggu_\infty$. Since $\s (\ggu_{\infty , n+1}) = \ggu_{\infty , n+1}$, the automorphism $\s$ induces the automorphism $\s_n$ of the Lie factor algebra $\ggu_{\infty } / \ggu_{\infty , n+1} \simeq \ggu_n$. Then $\s_n ( \pi_n ( I_{\l , n})) = \pi_n(I_\l )$, by Corollary \ref{c10Dec11}, where $\pi_n$ is as in (\ref{uinun}). Therefore, $\s ( I_{\l , n}) = I_{\l , n}$, as required.  $\Box $


Let $U_\infty := U(\ggu_\infty )$ be the universal enveloping algebra of the Lie algebra $\ggu_\infty$. The chain of Lie algebras $\ggu_2\subset \ggu_3\subset \cdots \subset \ggu_\infty = \cup_{n\geq 2} \ggu_n$ gives the chain of the universal enveloping algebras $U_2\subset U_3\subset \cdots \subset U_\infty = \cup_{n\geq 2} U_n$.

\begin{corollary}\label{b26Dec11}

\begin{enumerate}
\item  The inner derivations $\{ \ad (u)\, | \, u\in \ggu_\infty \}$ of the universal enveloping
algebra $U_\infty $ of the Lie algebra $\ggu_\infty $  are locally nilpotent
derivations. \item Every multiplicative subset $S$ of $U_\infty$ which is generated by an arbitrary set of elements of $\ggu_\infty$ is a
(left and right) Ore set in $U_\infty $. Therefore, $S^{-1}U_\infty \simeq
U_\infty S^{-1}$.
\end{enumerate}
\end{corollary}

{\it Proof}.  Both statements follow at once from Corollary \ref{b8Dec11} and the fact that $U_\infty = \cup_{n\geq 2} U_n$. $\Box $



\section{The  Lie algebra $\hggu_\infty$}\label{HHLIU}

In this section, the  completion $\hggu_\infty$ of the Lie algebra $\ggu_\infty$ with respect to the ideal topology on $\ggu_\infty$ is studied. Its properties diverge further from those of the Lie algebras $\ggu_n$ $(n\geq 2)$ and $\ggu_\infty$. For example, none of the nonzero inner derivations of the Lie algebra $\hggu_\infty$ is a locally nilpotent derivation; the Lie algebra $\hggu_\infty$ is neither locally nilpotent nor locally finite dimensional. The main result of this section is a classification of all the closed and open ideals of the topological Lie algebra $\hggu_\infty$ (Theorem \ref{20Jan12}.(1)). As a result, we proved that all the  open and all the closed ideals of the topological Lie  algebra $\hggu_\infty$ are topologically characteristic ideals (Corollary \ref{a21Jan12}.(3)).

The Lie algebra $\ggu_\infty$ is a topological Lie algebra where the topology is the {\em ideal topology}  on $\ggu_\infty$, that is, the set $\{ u+I\, | \, u\in \ggu_\infty ,  I\in \CJ(\ggu_\infty )\}$ is the basis of the ideal topology on $\ggu_\infty$ where $\CJ (\ggu_\infty )$ is the set of nonzero ideals of the Lie algebra $\ggu_\infty$.  Recall that the topological Lie algebra means that the maps $K \times \ggu_\infty \ra \ggu_\infty$, $ (\l ,v)\mapsto \l v$,
 $\ggu_\infty \times \ggu_\infty \ra \ggu_\infty$, $ (u,v)\mapsto u-v$, and $\ggu_\infty \times \ggu_\infty \ra \ggu_\infty$, $ (u,v)\mapsto [u,v]$, are continuous maps where the topologies on $K \times \ggu_\infty$ and $\ggu_\infty \times \ggu_\infty$ are the product topologies and the topology on $K$ is the discrete topology. By Theorem \ref{26Dec11}, every nonzero ideal of the Lie algebra $\ggu_\infty$ contains the ideal $\ggu_{\infty , n}$ for some $n\geq 1$. So, in the definition of the topology on $\ggu_\infty$ instead of all the nonzero ideals we can take the ideals $\{ \ggu_{\infty , n}\}_{n\geq 1}$.

The completion $\hggu_\infty$ of the topological space $\ggu_\infty$ is a topological Lie algebra
$$ \hggu_\infty =\{ \sum_{i\geq 1}a_i\der_i \, | \, a_i\in P_{i-1}, i\geq 1\}$$ where $ \sum_{i\geq 1}a_i\der_i$ is an infinite sum which is uniquely determined by its coefficients $a_i$.  The inclusion $\ggu_\infty \subseteq \hggu_\infty$ is an  inclusion of topological  Lie algebras, and the topology on $\hggu_\infty$ is (by definition) the strongest topology on $\hggu_\infty$ such that the map $\ggu_\infty \ra \hggu_\infty$, $a\mapsto a$, is continuous.

There is the strictly descending chain of ideals of the Lie algebra $\hggu_\infty$,
\begin{equation}\label{1infseru}
\hggu_{\infty , 1}= \hggu_{\infty} \supset \hggu_{\infty , 2} \supset \cdots \supset \hggu_{\infty , n}:=\{ \sum_{i\geq n} a_i\der_i \in \hggu_\infty \} \supset \cdots \supset \bigcap_{i\geq 1} \hggu_{\infty , i}=0
\end{equation}
and $\hggu_\infty / \hggu_{\infty , n+1}\simeq \ggu_\infty / \ggu_{\infty , n+1}\simeq  \ggu_{n+1} / P_n\der_{n+1}\simeq \ggu_n$ for all $n\geq 2$. The Lie algebra $\hggu_\infty$ can also be seen as the projective limit of the projective system of Lie algebra epimorphisms
$$ \cdots \ra \ggu_{n+1}\ra \ggu_n\ra \cdots \ra \ggu_3\ra \ggu_2\ra 0$$
where the epimorphism $\ggu_{n+1}\ra \ggu_n$ is the composition of the natural epimorphism $\ggu_{n+1}\ra \ggu_{n+1}/P_n\der_{n+1}$ and the isomorphism $\ggu_{n+1}/P_n\der_{n+1}\simeq \ggu_n$.
Each closed ideal $\hggu_{\infty , n}$ of the Lie algebra $\hggu_\infty$ is the completion/closure of the ideal $\ggu_{\infty ,n}$ of the Lie algebra $\ggu_\infty$. For each element $a\in \ggu_\infty$, the set $\{ a+\hggu_{\infty , n}\}_{n\geq 1}$ is the basis of all the open neighbourhoods of the element $a$. These sets $a+\hggu_{\infty , n}$, $n\geq 1$, are open and closed sets in $\hggu_\infty$.

The polynomial algebra $P_\infty =\cup_{n\geq 1}P_n=K[x_1, x_2, \ldots , ]$ is a left $\hggu_\infty$-module: for any $p\in P_n$ and $a=\sum_{i\geq 1} a_i\der_i \in \hggu_\infty$, $a\cdot p = \sum_{i=1}^n a_i\der_i (p)$.

\begin{proposition}\label{a19Jan12}

\begin{enumerate}
\item The Lie algebra $\hggu_\infty$ is not a solvable Lie algebra.
\item  The Lie algebra $\hggu_\infty$ is not a locally nilpotent  and not a locally finite dimensional Lie algebra.
\item Each element $u\in \hggu_\infty$ acts locally nilpotently on the $\hggu_\infty$-module $P_\infty$.
\item The chain of nonzero ideals in (\ref{1infseru}) is the derived  series for the Lie algebra $\hggu_\infty$, that is $(\hggu_\infty)_{(i)}=\hggu_{\infty , i+1}$ for all $i\geq 0$.
\item The upper central series for the Lie algebra $\hggu_\infty$ stabilizers at the first step, that is $(\hggu_{\infty})^{(0)}=\hggu_{\infty}$ and $(\hggu_{\infty})^{(i)}=\hggu_{\infty , 2}$ for all $i\geq 1$.
\item All the nonzero inner derivations of the Lie algebra $\hggu_\infty$ are not locally nilpotent derivations.
    \item $\Cen_{\hggu_\infty}(\der_1, \der_2, \ldots )=\{ \sum_{i\geq 1} \l_i\der_i\in \hggu_\infty\, | \, \l_i\in K\}$ is a maximal abelian Lie subalgebra of $\hggu_\infty$.
\item The centre $Z(\hggu_\infty )$ of the Lie algebra $\hggu_\infty$ is $0$.
\item  The Lie algebras  $\ggu_n$ where $n\geq 2$, $\ggu_\infty$ and $\hggu_\infty$  and are not pairwise isomorphic.
\end{enumerate}
\end{proposition}

{\it Proof}. 1. Statement 1 follows from statement 4.

2.  Let $a=\sum_{n\geq 1} \frac{x_1^n}{n!}\der_n$. Then
$$b_i:=(\ad \, \der_1)^i(a)=\sum_{n\geq i}\frac{x_1^{n-i}}{(n-i)!}\der_n\neq 0\;\; {\rm for \; all}\;\; i\geq 1.$$ So, the elements $b_1, b_2, \ldots $ are $K$-linearly independent. Therefore, the Lie subalgebra $\CG$ of $\hggu_\infty$ generated by the elements $\der_1$ and $a$ are neither nilpotent nor finite dimensional.

3. Statement 3 follows from Proposition \ref{25Dec11}.(3).

4.  Let $\d_i:=\ad (\der_i)$ for $i\geq 1$. Statement 4  follows from the following obvious facts: $\der_i\in \hggu_{\infty , i}$ and $[\der_i , \hggu_{\infty , i} ]=\hggu_{\infty , i+1}$ for all $i\geq 1$. These immediately imply that $(\hggu_\infty)_{(i)}\supseteq  \hggu_{\infty , i+1}$. The reverse inclusion follows from Proposition \ref{25Dec11}.(4) and the fact that $\hggu_\infty / \hggu_{\infty , i+1}\simeq \ggu_n$.

5.  Statement 5 follows from the facts that $[\der_1, \hggu_\infty ] = \hggu_{\infty , 2}$ and $[\der_1, \hggu_{\infty , 2}]=\hggu_{\infty , 2}$. The first equality yields the equality $[\hggu_\infty , \hggu_\infty ] = \hggu_{\infty , 2}$ and the second does the equalities  $(\hggu_\infty )^{(i)}=\hggu_{\infty , 2}$ for $i\geq 1$.

7.  $C:=\Cen_{\hggu_\infty}(\der_1, \der_2, \ldots )=\{ \sum_{i\geq 1} a_i\der_i\in \hggu_\infty \, | \, a_1\in K, \; a_i\in \cap_{j=1}^{i-1}\ker_{P_{i-1}}(\der_j)=K$ for all $i\geq 2\}$. Since the centralizer $C=\Cen_{\hggu_\infty}(C)$ is an abelian Lie subalgebra of $\hggu_\infty$, it is automatically a maximal abelian subalgebra of $\hggu_\infty$.

8.  Let $z=\sum_{i\geq 1}\l_i\der_i$ be a central element of the Lie algebra $\hggu_\infty$. By statement 7, $\l_i\in K$ for all $i\geq 1$. For all $i\geq 1$, $\l_i\der_{i+1} = [z, x_i\der_{i+1}]=0$,  i.e.,  $z=0$.

9. We know already that the Lie algebras $\ggu_n$ where $n\geq 2$ and $\ggu_\infty$ are not pairwise isomorphic (Proposition \ref{25Dec11}.(8)). The Lie algebras $\ggu_n$ are solvable (Proposition \ref{a8Dec11}.(1)) but the Lie algebra $\hggu_\infty$ is not (statement 1), hence $\hggu_\infty \not\simeq \ggu_n$ for all $n\geq 2$. The Lie algebra $\ggu_\infty$ is locally nilpotent (Proposition \ref{25Dec11}.(2)) but the Lie algebra $\hggu_\infty$ is not (statement 2), hence $\hggu_\infty \not\simeq \ggu_\infty$.

6.  Since the centre of the Lie algebra $\hggu_\infty$ is $0$ (statement 8) we have to show that if $a$ is a nonzero element of $\hggu_\infty$ then the inner derivation $\d =\ad (a)$ is not a locally  nilpotent  map.  Write $a=\sum_{i\in I}p_i\der_i= p_n\der_n+\cdots $ where  $p_i\in P_{i-1}$,  $I=\{ i\geq 1 \, | \, p_i\neq 0\} \neq \emptyset$ and $n$ is the minimal number in $I$. Let $b=\sum_{m\geq 0} \frac{x_n^m}{m!}\der_{n+m}$. Then
$$\d^i (b) = p_n^i\sum_{m\geq i} \frac{x_n^{m-i}}{(m-i)!}\der_{n+m}\neq 0, \;\; i\geq 1.$$
So, $\d$ is not a nilpotent map.  $\Box $


The next theorem gives a classification of all the open and the closed ideals of the topological Lie algebra $\hggu_\infty$, every nonzero closed ideal of $\hggu_\infty$ is open and vice versa (see statement 2 below).

\begin{theorem}\label{20Jan12}
\begin{enumerate}
\item  The set $\CJ (\hggu_\infty )$ of all the  open  ideals of the Lie algebra $\hggu_\infty$ is equal to the set  of all the nonzero closed ideals of the Lie algebra $\hggu_\infty$  and is equal to the set $\{ \hggu_\infty , \hggu_{\infty , n}, I_{\l , n}:= I_\l (n) +\hggu_{\infty , n+1}\, | \, n\geq 2, \l \in [1,\o^{n-1})\}$ where $I_\l ( n)$ is the ideal $I_\l$ of the Lie algebra $\ggu_n$ as defined in (\ref{Jlid}), that is $I_\l (n) = \oplus_{(\alpha , n)\leq \l} KX_{\alpha , n}$. In particular, every nonzero closed ideal $I$ of the Lie algebra $\hggu_\infty$ contains an ideal $\hggu_{\infty , n'+1}$ for some $n'\geq 2$, and, for $I=I_{\l , n}$, $n=\min \{ n'\geq 2\, | \, \hggu_{n'+1}\subseteq I\}$.
\item Let $I$ be a nonzero ideal of the Lie algebra $\hggu_\infty$. Then $I$ is an open ideal iff the ideal $I$ is a closed ideal.
\item The topological Lie algebras $\ggu_\infty$ and $\hggu_\infty$ are Hausdorff.
    \item The zero ideal is a closed but not open ideal of the topological Lie algebra $\hggu_\infty$.
\end{enumerate}
\end{theorem}

{\it Proof}. 1. We prove statement 1 for open ideals only. For the closed ones, it will follow from statement 2. Every open ideal of $\hggu_\infty$ is necessarily a nonzero ideal. It is obvious that the ideals $\hggu_\infty$, $\hggu_{\infty , n}$ and $I_{\l , n}$ in statement 1 are open ideals.
  Let $I$ be an open ideal of the Lie algebra $\hggu_\infty$ such that $I\neq \hggu_{\infty , n}$ for all $n\geq 1$. We have to show that $I=I_{\l, n}$ for some $\l$ and $n$ as in statement 1.
  The ideal $I$  contains an open neighbourhood  of the zero element, i.e.,  $\hggu_{\infty , n+1}\subseteq I$ for some natural element $n\geq 2$. We can assume that $n=\min \{ n'\, | \, n'\geq 2, \hggu_{\infty , n'+1}\subseteq I\}$. In view of the Lie algebra isomorphism, $\hggu_\infty / \hggu_{\infty , n+1}\simeq \ggu_n$, Theorem \ref{10Dec11}.(1), and the minimality of $n$, the claim and statement 1 are obvious.

  2. $(\Rightarrow )$ Let $I$ be an open ideal. By statement 1, which is true for open ideals, the ideal $I$ is as in statement 1.  All the ideals in statement 1, $\hggu_\infty $, $\hggu_{\infty , n}$ and $I_{\l , n}$, are obviously nonzero  closed ideals. Hence, $I$ is a nonzero closed ideal.

   $(\Leftarrow )$ Let $I$ be a nonzero closed ideal in $\hggu_\infty$. It suffices to show that $\hggu_{\infty , n+1}\subseteq I$ for some $n\geq 1$. This would automatically imply that $I$ is an open ideal. The ideal $\hggu_{\infty , n+1}$ is the closure of the ideal $\ggu_{\infty , n+1}$ of the Lie algebra $\ggu_\infty$. Fix a nonzero element of $I$, say $a=\sum_{m\geq n} a_m\der_m$ where $n\geq 1$, $a_m\in P_{m-1}$ for all $m\geq n$ and $a_n\neq 0$. Let $d$ be the total degree of the polynomial $a_n=\sum_{\alpha \in \N^{n-1}}\l_\alpha x^\alpha\in P_{n-1}$ where $\l_\alpha \in K$. Fix $\alpha = (\alpha_i) \in \N^{n-1}$ such that $|\alpha | := \alpha_1+\cdots +\alpha_{n-1}=d$. Applying $\ad (\der )^\alpha := \prod_{i=1}^{n-1}\ad (\der_i)^{\alpha_i}$ to the element $a$ yields the element of the ideal $I$ of the type $\alpha !\l_\alpha \der_n+\cdots $. So, without loss of generality we may assume from the very beginning that $a_n=1$, that is $a=\der_n+a_{n+1}\der_{n+1}+\cdots $.
 For each $\alpha \in \N^n$,
$$ I\ni  [ a,\frac{x^{\alpha +e_n}}{ \alpha_n+1}\der_{n+1}]=X_{\alpha , n +1}, $$
and so $P_n\der_{n+1}\subseteq I$.  For all natural numbers $l\geq n+1$, the ideal of the Lie algebra $\ggu_l$ generated by its subspace $P_n\der_{n+1}$ is equal to $\oplus_{i=n+1}^lP_{i-1}\der_i$. Therefore, $\ggu_{\infty , n+1} = \oplus_{i\geq n+1}P_{i-1}\der_i \subseteq I$. Hence, the closure of the ideal $\ggu_{\infty , n+1}$, which is $\hggu_{\infty , n+1}$, belong to $I$, as required.

  3. Let $a=\sum_{i\geq 1} a_i\der_i$ and $b=\sum_{i\geq 1} b_i\der_i$ be distinct elements of $\ggu_\infty$ (resp. $\hggu_\infty$) where $a_i, b_i\in P_{i-1}$ for all $i\geq 1$. Then $a_n\neq b_n$ for some $n$, hence $(a+\ggu_{\infty, n+1})\cap (b+\ggu_{\infty, n+1})=\emptyset$ (resp. $(a+\hggu_{\infty, n+1})\cap (b+\hggu_{\infty, n+1})=\emptyset$). This means that $\ggu_\infty$ and $\hggu_\infty$ are Hausdorff.

4. It is obvious.  $\Box$


A topological Lie algebra $\CG$ is called {\em closed uniserial} (resp. {\em open uniserial}) if the set of all closed (resp. open) ideals is a well-ordered set by inclusion. By Theorem \ref{20Jan12}.(1), the Lie algebra $\hggu_\infty$ is  a closed and open uniserial Lie algebra. So, the chain
\begin{equation}\label{1uinfid}
\hggu_\infty = \hggu_{\infty , 1}\supset \cdots \supset I_{\l , 2}\supset \cdots \supset I_{1,2}\supset \hggu_{\infty , 2}\supset \cdots \supset \hggu_{\infty, n}\supset \cdots \supset I_{\mu , n}\supset \cdots \supset I_{1, n}\supset \hggu_{\infty , n+1}\supset \cdots
\end{equation}
contains all the open and all the  nonzero closed ideals of the Lie algebra $\hggu_\infty$.


A topological Lie algebra $\CG$ is called an {\em open Artinian Lie algebra} (resp. a {\em closed Artinian Lie algebra}) if the set of open (resp. closed) ideals of $\CG$ satisfies the descending chain condition. A topological Lie algebra $\CG$ is called an {\em open Noetherian Lie algebra} (resp. a {\em closed Noetherian Lie algebra}) if the set of open (resp. closed) ideals of $\CG$ satisfies the ascending chain condition.  A topological Lie algebra $\CG$ is called an {\em open almost Artinian Lie algebra} (resp. a {\em closed almost Artinian Lie algebra}) if for all nonzero open (resp. closed) ideals $I$ of $\CG$ the Lie factor algebra $\CG / I$ is an Artinian Lie algebra. For a topological Lie algebra $\CG$, let 
$$\Aut_c(\CG):= \Aut_{{\rm Lie}}(\CG ) \cap \Aut_{\rm top}(\CG )$$ be its  group of automorphisms. Every element of $\Aut_c(\CG )$ is an isomorphism of the Lie algebra $\CG$ and an isomorphism of the topological space $\CG$. An ideal $I$ of a topological Lie algebra $\CG$ is called a {\em topologically characteristic ideal} if $\s (I)=I$ for all $\s\in \Aut_c(\CG )$.

 \begin{corollary}\label{a21Jan12}
\begin{enumerate}
\item The topological Lie algebra $\hggu_\infty$ is an open uniserial, closed uniserial, open almost Artinian and closed almost Artinian Lie algebra which is neither open nor closed Artinian and is neither open nor closed  Noetherian.
\item The uniserial dimensions of the sets of open and of closed ideals of the topological Lie algebra $\hggu_\infty$ coincide and are equal  to $\udim (\hggu_\infty ) = \o^\o$.
\item All the open/closed ideals of the Lie algebra $\hggu_\infty$ are topologically characteristic ideals.
\end{enumerate}
\end{corollary}

{\it Proof}. 1. Statement 1 follows from Theorem \ref{20Jan12}.(1).

2. Statement 2 follows from Theorem \ref{20Jan12}.(1) and Corollary \ref{a26Dec11}.(1).

3. The ideals $\hggu_{\infty , n}$ $(n\geq 1)$ are topologically characteristic ideals as they form the derived  series of the Lie algebra $\hggu_\infty$ (Proposition \ref{a19Jan12}.(4)). By Theorem \ref{20Jan12}.(1), it remains to show that every ideal $I_{\l , n}$ of $\hggu_\infty$ is a topologically characteristic ideal. Let $\s\in \Aut_c(\hggu_\infty )$. Since $\s (\hggu_{\infty , n+1})=\hggu_{\infty , n+1}$, the automorphism $\s$ induces the automorphism $\s$ of the Lie factor algebra $\hggu_\infty /\hggu_{\infty , n+1}\simeq \ggu_n$. Since $\hggu_{\infty , n+1}\subseteq I_{\l , n}$ and every ideal of the Lie algebra $\ggu_n$ is a characteristic ideal (Corollary \ref{c10Dec11}) we must have $\s (I_{\l , n}) = I_{\l , n}$, i.e.,  $I_{\l , n}$ is a topologically characteristic ideal of $\hggu_\infty$.  $\Box $



$${\bf Acknowledgements}$$

 The work is partly supported by  the Royal Society  and EPSRC.

\small{

Department of Pure Mathematics

University of Sheffield

Hicks Building

Sheffield S3 7RH

UK

email: v.bavula@sheffield.ac.uk}

\end{document}